\documentclass[11pt]{amsart}

\usepackage{amssymb,amsmath,amsthm}

\usepackage{amscd}
\newcommand{\grading}{\epsilon}
\newcommand{\tensor}{\mathbin{\widetilde{\otimes}}}
\newcommand{\semitensor}{\mathbin{\widehat{\otimes}}}
\newcommand{\conv}{\ast} % convolution symbol
\newcommand{\dist}{\mathcal{E}}

\newcommand{\dirac}{\mbox{$\not\negthinspace\partial$}}
\newcommand{\Z}{\mathbb{Z}}
\newcommand{\R}{\mathbb{R}}
\newcommand{\C}{\mathbb{C}}
\newcommand{\g}{\mathfrak{g}}
\newcommand{\h}{\mathfrak{h}}

\DeclareMathOperator{\sdim}{sdim}
\DeclareMathOperator{\ad}{ad}
\DeclareMathOperator{\Ad}{Ad}
\DeclareMathOperator{\Cl}{\mathbb{C}l}
\DeclareMathOperator{\ClR}{Cl}
\DeclareMathOperator{\Sym}{Sym}

\DeclareMathOperator{\End}{End}
\DeclareMathOperator{\Aut}{Aut}
\DeclareMathOperator{\Ext}{Ext}

\newcommand{\SR}[2]{\sideset{^{#1}}{^{#2}}{\mathop{\!SR}}}
\newcommand{\SRh}[2]{\sideset{^{#1}}{^{#2}}{\mathop{\!\widehat{SR}}}}

\newcommand{\K}[3]{\sideset{^{#1}}{^{\,#2}_{#3}}{\mathop{\!K\!}}(\mathrm{pt})}

\newcommand{\pair}[1][\,,]{\mathopen{\langle\!\langle\,} #1 \mathclose{\,\rangle\!\rangle}}

\DeclareMathOperator{\Spin}{Spin}
\DeclareMathOperator{\SO}{SO}
\DeclareMathOperator{\Hom}{Hom}
\DeclareMathOperator{\so}{\mathfrak{so}}
\DeclareMathOperator{\spin}{\mathfrak{spin}}
\DeclareMathOperator{\Spinc}{Spin^{c}}
\DeclareMathOperator{\U}{U}
\DeclareMathOperator{\SU}{SU}
\DeclareMathOperator{\PU}{PU}

\newtheorem{theorem}{Theorem}
\newtheorem{proposition}[theorem]{Proposition}
\newtheorem{lemma}[theorem]{Lemma}
\newtheorem*{schur}{Schur's Lemma}

\theoremstyle{remark}
\newtheorem*{definition}{Definition}
\newtheorem*{remark}{Remark}
\newtheorem*{example}{Example}

\title{Twisted Representation Rings and Dirac Induction}
\author{Gregory D. Landweber}

\email{greg@math.uoregon.edu}
\urladdr{http://www.uoregon.edu/\~{}greg/}
\address{Mathematics Department\\
         University of Oregon\\
         Eugene, OR 97403-1222}

\begin{document}

\begin{abstract}
Extending ideas of twisted equivariant $K$-theory, we construct twisted versions of the representation rings for Lie superalgebras and Lie supergroups, built from projective $\Z_{2}$-graded representations with a given cocycle.
We then investigate the pullback and pushforward maps on these representation rings (and their completions) associated to homomorphisms of Lie superalgebras and Lie supergroups.  As an application, we
consider the Lie supergroup $\Pi (T^{*}G)$, obtained by taking the cotangent
bundle of a compact Lie group and reversing the parity of its fibers.  An inclusion $H \hookrightarrow G$ induces a homomorphism
from the twisted representation ring of $\Pi(T^{*}H)$ to the twisted representation ring of $\Pi(T^{*}G)$, which pulls back via an algebraic version of the Thom isomorphism to give an additive homomorphism from $K_{H}(\mathrm{pt})$ to $K_{G}(\mathrm{pt})$ (possibly with twistings).  We then show that this homomorphism is in fact Dirac induction, which takes an $H$-module $U$ to the $G$-equivariant index of the Dirac operator $\dirac \otimes U$ on the homogeneous space $G/H$ with values in the homogeneous bundle induced by $U$.
\end{abstract}

\maketitle

%\tableofcontents

\section{Introduction}

Given a compact Lie group $G$ or a finite dimensional, reductive Lie algebra $\g$, the equivariant $K$-theory $K_{G}(\mathrm{pt})$ or $K_{\g}(\mathrm{pt})$
is the Grothendieck ring of isomorphism classes of finite dimensional
complex representations of $G$ or $\g$ respectively (see \cite{AS}, \cite{S0}).  Here we extend this to the supersymmetric or $\Z_{2}$-graded setting, constructing representation rings for Lie supergroups and Lie superalgebras.  We build our Grothendieck groups from supermodules, or $\Z_{2}$-graded representations, and in doing so we obtain an involution $\Pi$ on the isomorphism classes of supermodules, called ``parity reversal'', which interchanges the two homogeneous $\Z_{2}$-degree components.
%This gives us two natural ways to extend the representation
%ring of a compact Lie group or Lie algebra:  In the representation rings $R(\g)$ and $R(G)$,
%we identify the classes $[V] = [\Pi V]$ for supermodules $V$ (this convention appears in the superalgebra literature in for example \cite{BK0}), while in the
%super representation rings $SR(\g)$ and $SR(G)$ we impose the relation $[\Pi V] = -[V]$.
One natural way to extend the representation ring of a compact Lie group or Lie algebra is to define super representation rings
$SR(\g)$ and $SR(G)$ in which we impose the relation $[\Pi V] = -[V]$ for supermodules $V$.\footnote{ 
Another way is to define representation rings $R(\g)$ and $R(G)$ where we identify the classes of $V$ and $\Pi V$. This convention appears in the superalgebra literature in for example \cite{BK0}. These two sign conventions yield rings with similar additive structures but distinct products, as described by the author in \cite{L}. All our definitions and results in Sections \ref{subsection-representation-ring} to \ref{restriction-induction} hold for the representation rings $R(\g)$ and $R(G)$, but the theorems of
Sections \ref{Weyl-GKRS} and \ref{section-dirac-induction} do not.}
Furthermore, in analogy to complex $K$-theory, we can introduce Clifford-algebraic degree shifts, extending the super representation rings to $\Z_{2}$-graded rings $SR^{*}(\g)$ and $SR^{*}(G)$.
%The representation ring and the $\Z_{2}$-graded super representation ring are isomorphic as additive groups (although not canonically), but possess distinct product structures.  For some of the details, the reader is referred to \cite{L}, in which the author examines closely the $K$-theoretic properties of super representation rings for Lie superalgebras, working not only over the complex numbers, but also over the reals and other fields.

In this paper, we introduce a twisted version of the super representation ring, built from projective representations with a given cocycle.  This is closely related to $K$-theory with local coefficients developed by Donovan and Karoubi in \cite{DK}, and more recently to twisted equivariant $K$-theory used by Freed, Hopkins, and Teleman in \cite{F1,FHT1}.  The twistings we consider correspond to one-dimensional central extensions of our Lie superalgebras or Lie supergroups,
and so are classified by elements in the cohomology $H^{2}$.  Given homomorphisms of Lie superalgebras or Lie supergroups, we construct the
corresponding pullback or restriction maps of the (twisted) super representation rings, showing that $SR$ is a contravariant functor from Lie superalgebras and Lie supergroups to $\Z_{2}$-graded rings. Furthermore, we also obtain a pushforward or induction map in the spirit of Bott's paper \cite{B}, which is an additive group homomorphism acting on the completions $\SRh{}{}$ of the super representation rings with respect to a natural bilinear pairing on $SR$.
%(We note that we do \emph{not} obtain such a pushforward for the representation ring $R$.)

We apply this material to the case where $G$ is a compact (connected) Lie group with Lie subgroup $H$.  If we reverse the parity of the fibers of the cotangent bundles, we obtain Lie supergroups $\Pi(T^{*}G)$ and $\Pi(T^{*}H)$, whose underlying even Lie groups are $G$ and $H$, and whose
odd part is given by the coadjoint representations $\g^{*}$ and $\h^{*}$ respectively.  In Proposition \ref{thom-group}, we construct algebraic Thom isomorphisms between the twisted super representation rings
$$
\SR{\tau_{G}}{*}(G) \xrightarrow{\cong} \SR{b}{*+\dim G}\bigl(\Pi(T^{*}G)\bigr),
	\qquad
  \SR{\tau_{H}}{*}(H) \xrightarrow{\cong} \SR{b}{*+\dim H}\bigl(\Pi(T^{*}H)\bigr),
$$
where the twisting $b$ corresponds to a choice of $\Ad$-invariant inner
product on $\g$, and the twistings $\tau_{G}$ and $\tau_{H}$ correspond
to the projective cocycles of the spin representation $\mathbb{S}_{\g}$ and $\mathbb{S}_{\h}$
of the Clifford algebras $\Cl(\g)$ and $\Cl(\h)$, viewed as projective representations of $G$ and $H$ respectively.

The inclusion $i: H\hookrightarrow G$ of the Lie groups, along with the invariant inner product $b$ on $\g$, gives us an inclusion $j : \Pi(T^{*}H) \hookrightarrow \Pi(T^{*}G)$ of the corresponding Lie supergroups.  We then consider the restriction and induction maps associated to the inclusions $i$ and $j$.  These maps do \emph{not} commute with the Thom isomorphism, but their failure to commute is quite interesting.  For the restriction map,
we show in Theorem \ref{supergroup-restriction} that $j^{*}: \SR{b}{\dim G}\bigl(\Pi(T^{*}G)\bigr) \to \SR{b}{\dim H}\bigl(\Pi(T^{*}H)\bigr)$ pulls back via the Thom isomorphisms to the map
\begin{equation}\label{eq:GKRS-map1}
[V] \in \K{\tau_{G}}{}{G} \longmapsto [i^{*}V \otimes \mathbb{S}_{0}] -
[i^{*}V \otimes \mathbb{S}_{1}]
\in \K{\tau_{H}}{}{H},
\end{equation}
where $\mathbb{S} = \mathbb{S}_{0} \oplus\mathbb{S}_{1}$ is the $\Z_{2}$-graded spin representation associated to the orthogonal complement of $\h$ in $\g$, viewed here as a projective representation of $H$.  We observe that this map (\ref{eq:GKRS-map1}) is precisely the map considered by Gross, Kostant, Ramond, and Sternberg in \cite{GKRS}, which associates to each irreducible $G$-module a multiplet of irreducible $H$-modules,
as we recall in Theorem \ref{GKRS}.

On the other hand, we show in Theorem \ref{dirac-induction} that the induction map for the Lie supergroups gives an additive group homomorphism
$j_{*}: \SR{b}{\dim H}\bigl(\Pi(T^{*}H)\bigr) \to \SR{b}{\dim G}\bigl(\Pi(T^{*}G)\bigr)$, which pulls back via the Thom isomorphisms
to the Dirac induction map:
\begin{equation}\label{eq:dirac-induction1}
U \in \K{\tau_{H}-i^{*}\tau_{G}}{}{H} \longmapsto \mathrm{Index}_{G}\,\dirac^{G/H}_{U}\in \K{}{}{G}.
\end{equation}
When the coset space $G/H$ is a $\mathrm{Spin}^{c}$ manifold, a projective representation $U$ of $H$ with the appropriate cocycle gives rise to a Dirac operator
\begin{equation*}
\dirac^{G/H}_{U} : \Gamma \bigl( G \times_{H} ( \mathbb{S}_{0}\otimes U ) \bigr) \to
\Gamma \bigl( G \times_{H} ( \mathbb{S}_{1}\otimes U ) \bigr)
\end{equation*}
with coefficients in the homogeneous vector bundle $G \times_{H} U$ associated to $U$.  This Dirac operator is elliptic and homogeneous, i.e., it commutes with the $G$-actions on its domain and codomain, and so we can use Bott's formula from \cite{B} (which we recall in Theorem \ref{bott}) for its $G$-equivariant index. This Dirac induction map (\ref{eq:dirac-induction1}) is closely related to the holomorphic induction map given by the Borel-Weil-Bott theorem (in fact, it differs only by a $\rho$-shift,
which comes from tensoring with the canonical complex line bundle for $G/H$), as we discuss at the end of Section \ref{section-dirac-induction}.  In addition, this Dirac induction map, when extended to representation rings of positive energy representations over loop groups, as in \cite{L2} and \cite{Te}, becomes a vital tool for understanding the Freed-Hopkins-Teleman isomorphism between the Verlinde algebra and the twisted equivariant $K$-theory ${}^{\tau}\!K_{G}(G)$ in \cite{FHT}.

Our discussion begins with some preliminaries regarding associative superalgebras.  We then introduce the super representation ring in Section \ref{subsection-representation-ring}, reviewing the relevant results from our paper \cite{L}.  Our treatment here differs from that of \cite{L} in several ways. In this paper we consider only the complex case, which allows us to make some significant simplifications, and we develop the theory first for the super representation groups of associative superalgebras. In Section \ref{superalgebras-supergroups} we describe the super representation rings for Lie superalgebras $\g$ and Lie groups $G$ by passing to the universal enveloping algebra $U(\g)$ and convolution algebra $\dist(G)$ of distributions respectively.
For Lie supergroups, we consider supermodules carrying compatible actions of both the underlying even Lie groups $G_{0}$ and associated Lie superalgebras $\g$, by taking supermodules over a quotient of the semi-direct tensor product of $\dist(G_{0})$ and $U(\g)$.
In Section \ref{twisted}, we introduce twistings, considering projective representations of groups $G$ and Lie superalgebras $\g$ classified by $H^{3}(BG;\Z)$ and $H^{2}(\g)$ respectively.  Using our description in terms of associative superalgebras, we define twisted universal enveloping algebras as our basis for constructing the twisted super representation groups.
The remainder of the paper is concerned with proving the Thom isomorphism theorem
in both its Lie algebra and Lie group variants in Section \ref{thom}, defining the pullback and pushforward homomorphisms in Section \ref{restriction-induction}, and then measuring carefully their failure to commute in the case of the Lie supergroups of the form $\Pi(T^{*}G)$ in the final two sections.

%\smallskip
%\noindent
%\textbf{Acknowledgements.}
%The commutative diagrams in this article were typeset using the ``diagrams''  \TeX\ macro package by Paul Taylor.

\section{Associative superalgebras}

In this paper we work over the complex numbers unless otherwise noted.
A super vector space is a $\Z_{2}$-graded vector space $V = V_{0} \oplus V_{1}$.
Let $|v|$ denote the $\Z_{2}$-degree of a homogeneous element $v\in V$.
%, and define the grading involution by $\grading : v \mapsto (-1)^{|v|}\,v$.
%(Here, the $F$ stands for ``fermion number'', which comes from physics, where even elements with $\grading = +1$ are called ``bosons'' and odd elements with $\grading = -1$ are called ``fermions''.)
 A superalgebra is a $\Z_{2}$-graded algebra $A = A_{0} \oplus A_{1}$, which satisfies $|ab| = |a| + |b|$ for homogeneous elements $a,b\in A$. Two such elements are said to supercommute if $ab = (-1)^{|a||b|}ba$,
and we call the superalgebra supercommutative if any two of its elements supercommute.  In this section we consider unital associative superalgebras, with identity element in the even component.

If $V$ is a super vector space, then $\End(V)$ is a unital associative superalgebra.
The even component $\End(V)_{0}$ consists of all maps which %commute with $\grading$ and thus
preserve the grading on $V$, and the odd
component $\End(V)_{1}$ consists of all maps which
%anti-commute with $\grading$ and thus
interchange $V_{0}$ and $V_{1}$. In general,
a homomorphism between two super vector spaces is called even if it
%commutes with $\grading$, and so
preserves the grading (i.e., is
$\Z_{2}$-equivariant), or odd if it
%anti-commutes with $\grading$, and so
reverses the grading.  In this paper, whenever we refer simply to a
homomorphism (or isomorphism, endomorphism, etc.), we mean an \emph{even}
homomorphism unless otherwise noted.  A representation of a unital associative superalgebra $A$
on a super vector space $V$ is an (even) homomorphism $r : A \to \End(V)$,
and we call $V$ an $A$-supermodule.  We then have $|r(a)\,v| = |a| + |v|$
for homogeneous elements $a\in A$ and $v\in V$.

Given two unital associative superalgebras $A$ and $B$, their graded tensor product $A \tensor B$ has for its underlying vector space $A \otimes B$, and its multiplication is given by
$$ (a_{1} \tensor b_{1})\,(a_{2} \tensor b_{2} )
 := (-1)^{|b_{1}|\,|a_{2}|} (a_{1} a_{2}) \tensor (b_{1}b_{2}),$$
for homogeneous elements $a_{1},a_{2}\in A$ and $b_{1},b_{2}\in B$.
Both $A$ and $B$ inject into their graded tensor product $A\tensor B$ as $A\tensor 1$ and $1\tensor B$ respectively. The elements of $A$ and $B$ then supercommute with each other in $A \tensor B$. Given an $A$-supermodule $V$
and a $B$-supermodule $W$, their (exterior) tensor product $V \otimes W$ becomes an $A\tensor B$-supermodule, with action
$$(a\tensor b) (v \otimes w) = (-1)^{|b|\,|v|} a(v) \otimes b(w),$$
for homogeneous elements $a\in A$, $b\in B$, $v\in V$, $w\in W$.

Let $\ClR(n)$ and $\Cl(n)$ denote the real and complex Clifford algebras respectively, given by $n$ generators $\{e_{1},\ldots,e_{n}\}$
with relations
\begin{equation*}
e_{i}^{2} = -1 \qquad \text{and} \qquad e_{i} \cdot e_{j} = - e_{j} \cdot e_{i} \text{ for $i\neq j$}.
\end{equation*}
More generally, given a vector space $V$ with a symmetric bilinear form $b$, we define the Clifford algebra $\ClR(V,b)$ by
\begin{equation}\label{eq:clifford}
\ClR(V,b) := T^{*}(V) / \bigl(v \cdot w + w \cdot v = -2\,b(v,w) \text{ for $v,w\in V$}\bigr),
\end{equation}
where the tensor algebra $T^{*}(V)$ is the free unital associative algebra generated by $V$.
If the bilinear form is clear from the context, we write simply $\ClR(V)$.
Since $T^{*}(V)$ is $\Z$-graded and the ideal generated by the Clifford relations is contained entirely in even degrees, the Clifford algebra is a superalgebra.
We refer to $\Z_{2}$-graded representations of Clifford algebras as Clifford supermodules.
We recall that $\ClR(p) \tensor \ClR(q) \cong \ClR(p+q)$. So, the tensor product
of a $\ClR(p)$-supermodule with a $\ClR(q)$-supermodule is a $\ClR(p+q)$-supermodule. In general, we have
$$\ClR(V,b_{V}) \tensor \ClR(W,b_{W}) \cong \ClR( V \oplus W,b_{V}\oplus b_{W})$$
for vector spaces $V$ and $W$ with symmetric bilinear forms $b_{V}$ and $b_{W}$ respectively, and thus the tensor product of a
$\ClR(V)$-supermodule with a $\ClR(W)$-supermodule is a 
$\ClR(V \oplus W)$-supermodule.

\begin{schur}
Let $A$ be a collection of even and odd operators acting irreducibly on a super vector space $V = V_{0} \oplus V_{1}$ over an algebraically closed field. The only even $A$-equivariant endomorphisms of $V$ are scalar multiples of the identity, and there are two possibilities for the odd $A$-equivariant endomorphisms of $V$:
\begin{description}
\item[Type M]There are no nonzero odd endomorphisms of $V$.
\item[Type Q]The odd endomorphisms of $V$ are scalar multiples
of a parity reversing involution.
\end{description}
\end{schur}

Given a super vector space $V = V_{0}\oplus V_{1}$, its parity reversal $\Pi V$
has the same underlying vector space but the opposite $\Z_{2}$-grading.  In other words, we have $(\Pi V)_{0} = V_{1}$ and $(\Pi V)_{1} = V_{0}$.
%In terms of the grading operator, reversing the parity takes $\grading$ to $-\grading$.
Parity reversal allows us to view an odd homomorphism from $V$ to $W$ as an even homomorphism from $V$ to $\Pi W$ or from $\Pi V$ to $W$.  In particular, we have $\End(V)_{1} \cong \Hom(V,\Pi V)_{0}$ and $\End(V)_{0} \cong \Hom(V,\Pi V)_{1}$.
So, the odd endomorphisms described in Schur's Lemma can be viewed as even homomorphisms from $V$ to its parity reversal $\Pi V$.  We can likewise reverse the parity of an $A$-supermodule: if $V$ is an $A$-supermodule, then $\Pi V$ carries the same $A$-action on its underlying vector space, but has the opposite $\Z_{2}$-grading. As a consequence of Schur's Lemma, if $V$ is an irreducible $A$-supermodule which is isomorphic to its own parity reversal, $V\cong\Pi V$, then in fact there exists an even isomorphism $\alpha : V \to \Pi V$ such that $(\Pi \alpha) \circ \alpha = \mathrm{Id}$.

%\begin{remark}
%We note that Schur's Lemma as stated above requires that the field be algebraically closed.  For fields which are not algebraically closed, Schur's Lemma may fail, and we may indeed have  supermodules which are isomorphic to their parity reversals, but which do not admit a parity reversing involution. For example,
%when considering real Clifford supermodules, this gives rise to the $\Z_{2}$-torsion in $KO^{-n}(\mathrm{pt})$ for $n\equiv 1,2 \bmod 8$. See \cite{L} for a further discussion.
%\end{remark}

\section{Representation rings}
\label{subsection-representation-ring}

Given a unital associative superalgebra $A$, the (even) isomorphism classes of finite dimensional $A$-super\-modules form an abelian semi-group, and we can construct its corresponding Grothendieck group. To define the representation group and the super representation group,
we further consider the action of the parity reversal operator $\Pi$.

\begin{definition}
Let $F(A)$ be the free abelian group generated by (even) isomorphism classes of finite dimensional $A$-supermodules, and let $I(A)$ be the subgroup generated by the classes $[U] - [V] + [W]$ whenever there exists a short exact sequence
$$ 0 \to U \to V \to W \to 0.$$
%The \emph{representation group} of $A$ is 
%$$R(A) := F(A) / I_{-}(A),$$
%where $I_{-}(A)$ is the subgroup generated by both $I(A)$ and anti-dual classes of the form $[V] - [\Pi V]$.\\
The \emph{super representation group} of $A$ is
$$SR(A) := F(A) / I_{+}(A),$$
where $I_{+}(A)$ is the subgroup generated by both $I(A)$ and self-dual classes $[V]$ for supermodules $V$ isomorphic to their own parity reversals, $V\cong \Pi V$.\footnote{
Alternatively, the \emph{representation group} of $A$, as discussed in \cite{BK0}, uses the opposite sign convention,
defining $R(A) := F(A) / I_{-}(A)$, where $I_{-}(A)$ is generated by both $I(A)$ and anti-dual classes of the form $[V] - [\Pi V]$.

Over fields for which Schur's Lemma fails, we define the super representation group more precisely as the quotient of $F(A)$ by the subgroup generated by $I(A)$ and classes $[V]$ for $V$ admitting parity reversing involutions. Without Schur's Lemma, a supermodule may be isomorphic to its parity reversal but not via an involution, giving rise to a 2-torsion subgroup measuring the failure of Schur's Lemma.  See \cite{L} for further discussion of this general case.}
\end{definition}

%We call this the ``representation group'' and ``super representation group'' in analogy to the terms ``trace'' and ``supertrace''.
The parity reversal operator $\Pi$ descends to the super representation group, with $[\Pi V] = - [V]$.
%We observe that the parity reversal operator $\Pi$ descends to the super representation group, acting by $-1$ on $SR(A)$.
%and it acts by $+1$ on $R(A)$ and by $-1$ on $SR(A)$. So, in the representation group, we identify the class of an $A$-module $V$ with that of its parity reversal $\Pi V$. (See \cite{BK0} for a discussion of this representation group.) However, in the super representation group, we have $[\Pi V] = - [V]$.
%The representation group and super representation group can be viewed as generalizations of the equivariant cohomology and equivariant $K$-theory of a point respectively (\cite{V}).
If $A = A_{0}$ is a purely even algebra, then $SR(A_{0}) \cong \K{}{}{A_{0}}$, where $\K{}{}{A_{0}}$ is the conventional representation group, the Grothendieck group of isomorphism classes of finite dimensional ungraded $A_{0}$-modules.  It is instructive to examine this isomorphism more closely.  Since $A_{0}$ has no odd component, an $A_{0}$-supermodule is just an (ordered) pair of ungraded $A_{0}$-modules $(V_{0},V_{1})$.
%In the representation group $R(A_{0})$,
%we symmetrize, identifying the classes $[(V_{0},V_{1})]$ and $[(V_{1},V_{0})]$. Such a pair then maps to the direct sum of its two components in $\K{}{}{A_{0}}$:
%\begin{equation}
%[(V_{0},V_{1})] \in R(A_{0}) \longmapsto [V_{0}\oplus V_{1}] = [V_{0}] + [V_{1}] \in \K{}{}{A_{0}}.
%\end{equation}
%On the other hand, 
The super representation group is anti-symmetric, with $[(V_{0},V_{1})] = -[(V_{1},V_{0})]$.  We can view the class of such a pair in $SR(A_{0})$ as a formal direct difference, or virtual $A_{0}$-module:
\begin{equation}\label{eq:direct-difference}
[(V_{0},V_{1})] \in SR(A_{0}) \mapsto \text{``}[V_{0}\ominus V_{1}]\text{''} = [V_{0}] - [V_{1}] \in \K{}{}{A_{0}}.
\end{equation}
In this regard, the super representation group closely resembles the Grothendieck construction.

%\begin{remark}
%When working over fields for which Schur's Lemma fails, we define the super representation group more precisely as the quotient of $F(A)$ by the subgroup generated by $I(A)$ and classes $[V]$ for $V$ admitting parity reversing involutions. Without Schur's Lemma, we may have supermodules which are isomorphic to their parity reversals but not via an involution, giving rise to a 2-torsion subgroup measuring the failure of Schur's Lemma.  See \cite{L} for further discussion of this general case.
%\end{remark}

%\begin{remark}
%We note that Schur's Lemma as stated above requires that the field be algebraically closed.  For fields which are not algebraically closed, Schur's Lemma may fail, and we may indeed have  supermodules which are isomorphic to their parity reversals, but which do not admit a parity reversing involution. For example,
%when considering real Clifford supermodules, this gives rise to the $\Z_{2}$-torsion in $KO^{-n}(\mathrm{pt})$ for $n\equiv 1,2 \bmod 8$. See \cite{L} for a further discussion.
%\end{remark}

We can introduce degree shifts into the super representation group by incorporating Clifford algebras via graded tensor products:

\begin{definition}
The $n$-times degree-shifted super representation group is:
$$SR^{-n}(A) := SR\bigl( A \tensor \Cl(n) \bigr),$$
constructed from $A$-supermodules admitting $n$ supercommuting supersymmetries.
\end{definition}

This Clifford-algebraic definition is motivated by \cite{ABS}, where Atiyah, Bott, and Shapiro establish the connection between Clifford algebras and $K$-theory, proving for the trivial algebra $A = \C$ that
\begin{equation}\label{eq:ABS}
	SR^{-n}(\C) = SR\bigl( \Cl(n) \bigr)
	\cong \widetilde{K}(S^{n}) = K^{-n}(\mathrm{pt}).
\end{equation}
%(Actually, they use different notation and slightly different definitions, as well as prove the analogous isomorphism for the real super representation ring and the $KO$-theory of a point.)
The twofold periodicity of complex Clifford algebras gives rise to a twofold periodicity of the degree-shifted super representation group by Morita equivalence, as we prove in \cite[\S 6]{L}:
%On the other hand, degree shifting does not affect the representation group at all!

\begin{proposition}
%The representation group is independent of the degree shift, with $R^{-n}(A) \cong R(A)$, while
The super representation groups have twofold periodicity: $SR^{-n}(A) \cong SR^{-n-2}(A)$.
%Furthermore, we have a non-canonical isomorphism $R(A) \cong SR^{0}(A) \oplus SR^{1}(A)$ as additive groups.
\end{proposition}

%These isomorphisms can be made explicit by considering canonical bases for the degree-shifted representation group and super representation group.  Since we are dealing with finite dimensional representations, the class of an $A$-supermodule in $R(A)$ or $SR(A)$ decomposes as a sum of classes of irreducible $A$-modules.  Both $R(A)$ and $SR(A)$ can then be described as free abelian groups on bases of irreducibles.  Recalling Schur's Lemma, the irreducibles come in two flavors: type M with no odd endomorphims, and type Q which admit an odd involution.  In terms of parity reversal, irreducibles $M$ of type M come in $M, \Pi M$ pairs, while irreducibles $Q$ of type Q satisfy $Q \cong \Pi Q$.  When we shift degrees, the roles are reversed. If $M$ is a degree 0 irreducible of type M, then $M \oplus \Pi M$ is a degree 1 irreducible of type Q. On the other hand, if $Q$ is a degree 0 irreducible of type Q, then $Q$ admits two distinct $\Cl(1)$ actions which are parity reversals of each other, giving a pair of degree 1 irreducibles $Q_{+}, Q_{-}$ of type M.
%For a basis of the representation group $R^{0}(A)$, we can take one element of each $[M],[\Pi M]$ pair and all $[Q]$, while for $R^{1}(A)$ we can take all $[M\oplus\Pi M]$ and one element from each $[Q_{+}], [Q_{-}]$ pair. These two groups are clearly isomorphic, with an isomorphism given on the canonical bases by
%\begin{align*}
%	[M] = [\Pi M]\in R^{0}(A) &\longmapsto [M \oplus \Pi M]\in R^{1}(A), \\
%	[Q] \in R^{0}(A) &\longmapsto [Q_{+}] = [Q_{-}]\in R^{1}(A).
%\end{align*}
The two components of the degree-shifted super representation group can be made explicit by constructing their (almost) canonical bases.  Since we are dealing with finite dimensional representations, the class of an $A$-supermodule in $SR(A)$ decomposes as a sum of classes of irreducible $A$-modules.  Then $SR(A)$ can be described as a free abelian group on a basis of irreducibles.  Recalling Schur's Lemma, the irreducibles come in two flavors: type M with no odd endomorphims, and type Q which admit an odd involution.  In terms of parity reversal, irreducibles $M$ of type M come in $M, \Pi M$ pairs, while irreducibles $Q$ of type Q satisfy $Q \cong \Pi Q$.  When we shift degrees, the roles are reversed. If $M$ is a degree 0 irreducible of type M, then $M \oplus \Pi M$ is a degree 1 irreducible of type Q. On the other hand, if $Q$ is a degree 0 irreducible of type Q, then $Q$ admits two distinct $\Cl(1)$ actions which are parity reversals of each other, giving a pair of degree 1 irreducibles $Q_{+}, Q_{-}$ of type M.
The super representation group is generated by only the type M irreducibles, so
a basis for $SR^{0}(A)$ is given by choosing one element from each $[M],[\Pi M]$
pair, while a basis for $SR^{1}(A)$ is given by choosing one element from each $[Q_{+}], [Q_{-}]$ pair.\footnote{
In contrast, the representation group $R(A)$ is unchanged under degree shifts, and it has a canonical basis consisting of the classes $[M] = [\Pi M]$ and $[Q]$. This gives us a non-canonical additive isomorphism between $R(A)$ and $SR^{*}(A)$. }
%The full $\Z_{2}$-graded super representation group $SR^{*}(A) = SR^{0}(A) \oplus SR^{1}(A)$ is therefore isomorphic to the representation group $R(A)$, with isomorphism given on the basis by
%\begin{align*}
%	[M] = -[\Pi M]\in SR^{0}(A) &\longmapsto [M] \in R(A), \\
%	[Q_{+}] = -[Q_{-}] \in SR^{1}(A) &\longmapsto [Q] \in R(A).
%\end{align*}
%Note that this isomorphism depends for its signs on the choice of preferred element in each $[M],[\Pi M]$ or $[Q_{+}],[Q_{-}]$ pair, and thus it is not canonical.

We note that if $A = A_{0}$ has no odd component, then the graded and ungraded tensor products with the Clifford algebra agree: $A_{0} \tensor \Cl(n) \cong A_{0} \otimes \Cl(n)$. In this
case, it follows that modules over $A_{0} \tensor \Cl(n)$ decompose as direct sums of modules of the form $V \otimes \mathbb{S}$, where $V$ is an ungraded $A_{0}$-module, and $\mathbb{S}$ is a $\Cl(n)$-supermodule.  In terms of representation rings, we find that
\begin{equation}\label{eq:conventional-bott}
SR^{-n}(A_{0}) \cong \K{}{}{A_{0}} \otimes SR^{-n}(\C)
\cong \K{}{}{A_{0}} \otimes K^{-n}(\mathrm{pt}) \cong
\begin{cases}
	\K{}{}{A_{0}} & \text{ for $n$ even}, \\
	0         & \text{ for $n$ odd},
\end{cases}
\end{equation}
recalling the Atiyah-Bott-Shapiro isomorphism (\ref{eq:ABS}).

If $A$ is a Hopf superalgebra, then its comultiplication homomorphism $\Delta: A \to A\tensor A$ gives us a ring structure on the super representation group as follows:
Given two $A$-supermodules $V$ and $W$,
their tensor product $V \otimes W$ is a representation of the graded tensor product $A \tensor A$.  Pulling back to $A$ via $\Delta$, we obtain an $A$-supermodule $\Delta^{*}(V\otimes W)$ which we call the interior tensor product of $V$ and $W$.  If the diagonal
map is clear from the context, then we write simply $V \otimes W$.  For the interior tensor product, the parity reversal operator obeys
$$\Pi (V \otimes W) \cong (\Pi V) \otimes W \cong V \otimes (\Pi W).$$
As a consequence, the subgroup $I_{+}(A)$ in the definition
of $SR(A)$ is an ideal in $F(A)$, and thus the 
super representation group is in fact a ring. The identity element in this ring is the class $I = [\C]$ of the trivial, purely even representation.  In terms of the Hopf superalgebra structure, the trivial representation is the counit homomorphism $A \to \C$, and the axioms for a Hopf superalgebra ensure that this product does indeed give a ring structure.\footnote{
In contrast, the preferred ring structure for the representation ring $R(A)$, as described for example in \cite{BK0}, is the tensor product, except that the product of two irreducibles of type Q is taken to be \emph{half} their tensor product. }

%\begin{remark}
%Actually, the preferred ring structure on the representation group $R(A)$, as described for example in \cite{BK0}, is not the
%one induced by the tensor product. In terms of the canonical basis described above, the product of two irreducibles of type M or of irreducibles of types M and Q
%is given by the standard tensor product. However, the class of the tensor product of two irreducibles of type Q is divisible by two in $R(A)$, and one defines their product in the representation group to be \emph{half} their interior tensor product.
%\end{remark}

If we start with a Hopf superalgebra $A$, when we introduce degree shifts, the graded tensor products $A \tensor \Cl(n)$ are not in general Hopf superalgebras, and so the degree-shifted super representation groups $SR^{-n}(A)$ do not necessarily carry individual ring structures.  However, if $V$ is an $A \tensor\Cl(n)$-supermodule
and $W$ is an $A\tensor\Cl(m)$-supermodule, then $V \otimes W$ is a supermodule over
$$\bigl(A \tensor \Cl(n)\bigr) \tensor \bigl(A \tensor\Cl(m)\bigr) 
\cong (A \tensor A) \tensor \Cl(n+m).$$
Restricting to the diagonal, the interior tensor product then gives an $A\tensor\Cl(n+m)$-supermodule, which furthermore satisfies the property
$V \otimes W \cong \Pi^{nm} (W \otimes V).$
This induces a product
$$SR^{-n}(A) \otimes SR^{-m}(A) \to SR^{-n-m}(A),$$ with respect
to which the full super representation ring $SR^{*}(A)$ becomes a supercommutative $\Z_{2}$-graded ring over the degree zero component
$SR^{0}(A) = SR(A)$ (see \cite[\S 6.2]{L} for a full discussion).

\section{Lie superalgebras and Lie supergroups}
\label{superalgebras-supergroups}

\subsection{Lie superalgebras}
A Lie superalgebra, as defined by Kac in \cite{Kac}, is a super vector space $\g =\g_{0} \oplus \g_{1}$ with a bilinear $\Z_{2}$-graded product $[\,,]\colon \g\times\g\to\g$, referred to as a bracket, satisfying
\begin{itemize}
\item $[X,Y] = - (-1)^{|X|\,|Y|}[Y,X]$,
\item $\bigl[X, [Y,Z] \bigr] = \bigl[[X,Y],Z\bigr]
		+ (-1)^{|X|\,|Y|}\bigl[Y,[X,Z]\bigr]$,
\end{itemize}
for homogeneous elements $X, Y, Z\in \g$. In other words, the bracket is super anti-commutative, and the adjoint action $\ad_{X} : Y \mapsto [X,Y]$ is a super derivation. Alternatively, a Lie superalgebra consists of a conventional even Lie algebra
$\g_{0}$ and an odd $\g_{0}$-module $\g_{1}$, equipped with a $\g_{0}$-invariant symmetric bilinear form $\g_{1}\otimes\g_{1}\to\g_{0}$.  An associative superalgebra $A$ always admits a Lie superalgebra structure by defining the
bracket $[a,b] := ab - (-1)^{|a|\,|b|}\,ba$ as the super commutator for homogeneous elements $a,b\in A$.

A Lie superalgebra is neither associative nor unital. However, given a Lie superalgebra $\g$, we can construct its universal enveloping algebra
$$U(\g) := T^{*}(\g) / \bigl(X Y - (-1)^{|X|\,|Y|}YX = [X,Y]\text{ for homogeneous $X,Y\in\g$} \bigr),$$
where the tensor algebra $T^{*}(\g)$ is the free unital associative algebra
generated by $\g$. The universal enveloping algebra is then a unital associative algebra which comes equipped with a canonical injection $\g \hookrightarrow U(\g)$, and we view $\g$ as a Lie sub-superalgebra of $U(\g)$.
Any Lie superalgebra homomorphism $\g \to A$, where $A$ is an associative algebra, factors uniquely through $U(\g)$, and in fact $U(\g)$ can be defined
via this universal property.
In particular, every representation $\g \to \End(V)$ lifts to an algebra homomorphism $U(\g) \to \End(V)$, and conversely every algebra homomorphism $U(\g)\to\End(V)$ restricts to a representation $\g \to \End(V)$.  Furthermore,
the univeral enveloping algebra $U(\g)$ is a Hopf superalgebra (see \cite{MM}).
In particular, the comultiplication homomorphism $\Delta : U(\g) \to U(\g) \tensor U(\g)$ is given by
$$\Delta X = X \tensor 1 + 1 \tensor X,$$
on the generators $X\in \g$.  This gives a Lie algebra homomorphism $\Delta : \g \to U(\g) \tensor U(\g)$, which then extends to $U(\g)$ by the universal property. The counit homomorphism $U(\g) \to \C$ is given on the generators
by $X \mapsto 0$, and it likewise extends to $U(\g)$.

\begin{definition}
The super representation ring of a Lie superalgebra $\g$ is the graded ring with components
	$$%R(\g) := R\bigl(U(\g)\bigr), \qquad
	SR^{-n}(\g) := SR^{-n}\bigl(U(\g)\bigr)
	= SR\bigl( U(\g) \tensor \Cl(n) \bigl),$$
constructed using representations of the universal enveloping algebra $U(\g)$.
\end{definition}

When considering degree shifts, each homogeneous component $SR^{-n}(\g)$ of the super representation ring is constructed from $\g$-supermodules carrying auxiliary $\Cl(n)$-actions. In addition, the action of the Clifford generators must commute with the action of the even component $\g_{0}$ and anti-commute with the action of the odd component $\g_{1}$ of the Lie superalgebra. We call such representation spaces $\g$-Clifford supermodules.

Here we are considering complex Lie superalgebras and their complex universal enveloping algebras.  In some of the sections that follow,
we work instead with \emph{real} Lie superalgebras, but ultimately we still
want to construct super representation rings in terms of supermodules over
a complex associative algebra. In what follows, if $\g$ is a real Lie superalgebra, then when we take its universal enveloping algebra, we
actually mean the complexification $U(\g\otimes\C) \cong U(\g) \otimes \C$.
 
\subsection{Lie groups}
\label{liegroups}
We can also construct super representation rings of Lie groups in terms of representations of unital associative algebras. Recall that any representation $r : G \to \Aut(V)$ of a \emph{finite} group $G$ lifts to a representation $r: \C[G] \to \End(V)$  of its complex group algebra by the formula
$$r: {\bigoplus}_{g \in G}{a_{g}\,[g]} \longmapsto {\sum}_{g\in G} a_{g}\,r(g),$$
for complex coefficients $a_{g} \in \C$ for each $g\in G$. Conversely, given a representation of the complex group algebra, we can reconstruct the representation of the underlying discrete group by taking
$r(g) = r([g])$, i.e., by setting the coefficients $a_{g} = 1$ and $a_{h} = 0$ for all $h\neq g$.
Like the universal enveloping algebra, the complex group algebra $\C[G]$ is also a Hopf algebra, with comultiplication homomorphism $[g] \mapsto [g] \otimes [g]$ induced by the diagonal map $\Delta : G \to G \times G$, and counit homomorphism $\bigoplus_{g\in G}a_{g}\,[g] \mapsto \sum_{g\in G}a_{g},$ summing all the coefficients.

If $G$ is a compact Lie group, then a representation is a continuous homomorphism $r : G \to \Aut(V)$, and instead of lifting to the group algebra, we consider the ring  $C^{\infty}(G)$ of smooth complex-valued functions on $G$.  The action of a function $f\in C^{\infty}(G)$ is then given by the integral
\begin{equation}\label{eq:cg-lift}
	r : f\in C^{\infty}(G) \longmapsto \int_{G} f(g)\,r(g)\,dg \in \End(V),
\end{equation}
where $dg$ is the bi-invariant Haar measure on $G$. The product we use on $C^{\infty}(G)$ is not simply multiplication of functions, but rather the convolution given by
\begin{equation}\label{eq:cg-product}
(f_{1}\conv f_{2})(g) := \int_{G}f_{1}(gh^{-1})\,f_{2}(h)\,dh = \int_{G} f_{1}(gh)\,f_{2}(h^{-1})\,dh
\end{equation}
for $f_{1},f_{2}\in C^{\infty}(G)$, which can be viewed as a generalization of the product on the group algebra.  We note that $C^{\infty}(G)$ need not have a unit with respect to this convolution product.

We actually work not directly with $C^{\infty}(G)$, but rather with a completion $\dist(G)$. Using the product (\ref{eq:cg-product}), any smooth function $f\in C^{\infty}(G)$ gives a linear convolution operator $f\conv{\bullet}: C^{\infty}(G) \to C^{\infty}(G)$ with kernel $f$. For our completion, we take $\dist(G)$ to consist of distributions on $G$ whose corresponding convolutions give operators in $\End( C^{\infty}(G))$ (see \cite{AM}). A representation $r: G \to \Aut(V)$ on a finite dimensional vector space $V$ is in fact a smooth group homomorphism, so we can view $r$ as an element of $C^{\infty}(G) \otimes \End(V)$.  The action of a distribution $\phi \in \dist(G)$ is still given by (\ref{eq:cg-lift}), or more precisely we take $r(\phi) = (\phi \conv r^{-1} )(e) \in \End(V)$. This gives a lift $r : \dist(G) \to \End(V)$ of our representation which is indeed an algebra homomorphism with respect to convolution.

Conversely, given a representation $r : \dist(G) \to \End(V)$, we can recover the underlying representation of $G$ by considering Dirac delta distributions $\delta_{g} \in \dist(G)$ for $g\in G$, defined by the identity
$$\int_{G}\delta_{g}(h)\,f(h)\,dh = f(g)$$
for all $f\in C^{\infty}(G)$. In terms of the convolution product, we have $\delta_{g} \conv f = l_{g}f$ for $f\in C^{\infty}(G)$, where the left $G$-action on $C^{\infty}(G)$ is given by $(l_{g}f)(h) = f(g^{-1}h)$. It follows that $\delta_{g} \conv \delta_{h} = \delta_{gh}$, and so we obtain an injective group homomorphism $\delta : G \hookrightarrow \dist(G)$.  Pulling back the representation $r$ we
obtain the original group representation as $r(g) = r(\delta_{g})$.

The completion $\dist(G)$ is a unital associative algebra with respect
to the convolution product, and its identity element is the Dirac delta distribution $\delta_{e}$ at the identity.  It is in fact a Hopf algebra,
with comultiplication induced by the diagonal map $\Delta : G \to G \times G$ and counit given by $\phi \mapsto (\phi \conv 1)(e)$, where $1$ here refers to the constant function with value $1$. In other words, the counit homomorphism maps a distribution $\phi$ to its average over the group $\int_{G}\phi(g)\,dg$.  This Hopf algebra structure allows us to define:

\begin{definition}The %representation ring and
super representation ring of a compact Lie group $G$ is the graded ring given by
\begin{align*}
%R(G) := R \bigl( \dist(G) \bigr), \qquad
SR^{-n}(G) := SR^{-n}\bigl( \dist(G) \bigr) = SR \bigl( \dist(G) \tensor \Cl(n) \bigr),
\end{align*}
constructed using representations of the convolution algebra of distributions $\dist(G)$.
\end{definition}

In other words, we consider representations of $G$ carrying auxiliary $\Cl(n)$-actions commuting with the action of $G$, which we call $G$-Clifford supermodules. We note that if $G$ is a compact Lie group (or $\g =\g_{0}$ is a purely even Lie algebra), then $\dist(G)$ (respectively $U(\g_{0})$) is purely even, and the Clifford component completely decouples from the Lie component as in (\ref{eq:conventional-bott}). So, we observe that this super representation ring construction is most interesting when applied in the supersymmetric setting, involving representations of Lie superalgebras, or as we shall discuss below, Lie supergroups.

\begin{remark}
Alternatively, we could reprise the above constructions using only representations of groups, rather than superalgebras.  To introduce degree shifts, we could consider $\Z_{2}$-graded representations of the group
$G \times \Cl(n)^{*}$, where $\Cl(n)^{*}$ is the group of invertible elements in the Clifford algebra.  Since $\Cl(n)^{*}$ is dense in $\Cl(n)$, its representations and those of
the Clifford algebra are interchangeable, as are representations of the Lie group $G$ and the unital associative algebra $\dist(G)$.
\end{remark}

\subsection{Lie supergroups}
\label{supergroups}

A supermanifold can be thought of as an underlying
even conventional manifold with odd ``fuzz'' (see \cite{F}
for an excellent exposition or \cite{W} for a thorough treatment of the subject),
and is best described in terms of
its ring of functions.  The functions on a supermanifold
form a $\Z_{2}$-graded supercommutative ring, where both
the even and odd components are modules over the ring of functions on the underlying even conventional manifold.
Locally, supermanifolds look like regions of $\R^{p|q} = \R^{p} \times \Pi \R^{q}$, where $\R^{p}$ is the
underlying even conventional manifold, and $\Pi \R^{q}$ is
the odd ``fuzz''.  The
functions on the even part $\R^{n}$ lie in the closure of the polynomial algebra $\Sym^{*}(\R^{n})$. On the other hand,
the functions on the odd part $\Pi \R^{m}$ are anti-commuting polynomials, which instead lie in the exterior
algebra $\Lambda^{*}(\R^{m}) = \Sym^{*}(\Pi \R^{n})$.
Since the exterior algebra is finite dimensional, there
is no need to consider completions. The ring of
functions on the supermanifold $\R^{p|q}$ is therefore $C^{\infty}(\R^{p|q}) = C^{\infty}(\R^{p})\otimes\Lambda^{*}(\R^{q})$, with even part $C^{\infty}(\R^{p})\otimes\Lambda^{\text{even}}(\R^{q})$ and odd
part $C^{\infty }(\R^{p})\otimes\Lambda^{\text{odd}}(\R^{q})$.

An example of a supermanifold is $\Pi(TX)$, the
parity reversed tangent bundle of a smooth manifold $X$,
where the fibers are treated as odd vector spaces.
Its ring of smooth functions is
$$C^{\infty}\bigl(\Pi(TX)\bigr) = \Gamma\bigl(\Lambda^{*}(T^{*}X)\bigr) = \Omega^{*}(X),$$
the space of differential forms on $X$. In general, if $E\to X$ is a smooth real vector
bundle, then $\Pi E$ is a supermanifold with the ring of functions
$C^{\infty}(\Pi E) = \Gamma\bigl(\Lambda^{*}(E^{*})\bigr)$, i.e., sections of the bundle whose fibers are the exterior
algebras of the fibers of the dual bundle $E^{*}$, or in other words, anti-commuting polynomials on the fibers of $E$ whose coefficients are functions on the base $X$.

A Lie supergroup is a supermanifold which carries a smooth $\Z_{2}$-graded group structure (see \cite{W} or \cite{DM} for a complete definition). The underlying even conventional manifold of a Lie supergroup is a Lie group, and its odd component is a vector space with no further topological or geometric structure.
The space of left-invariant vector fields on a Lie supergroup is a Lie superalgebra.
If $G$ is a Lie supergroup, then the inclusion of the identity $e\hookrightarrow G$ induces a ring homomorphism $C^{\infty}(G)\to C^{\infty}(e)\cong\R$, and we denote this evaluation map by $f\mapsto f(e)$.
The tangent space at the identity $TG_{e}$ is the vector
space of all super-derivations $D\colon C^{\infty}(G)\to C^{\infty}(e)\cong\R$ satisfying 
\begin{equation*}
	D(fg) = D(f)\,g(e) + (-1)^{|D|\,|f|} f(e)\,D(g).
\end{equation*}
Note that these super-derivations include even derivations corresponding to tangent vectors of the underlying even Lie group, as well as odd derivations corresponding to vectors
in the directions of the odd ``fuzz''. This tangent space $TG_{e}$ is then a
Lie superalgebra.

If $G_{0}$ is a Lie group and $V$ is a real representation
of $G_{0}$, then we can consider the vector bundle $E = (G_{0}\times G_{0}) \times_{\Delta G_{0}} V$ over $G_{0}$ with fibers $V$.
Then $\Pi E$ is a Lie supergroup with Lie superalgebra $\g_{0} \oplus \Pi V$.
%(Such bundles can be trivialized by left translation as $E \cong G_{0}\times V$.
%However, when we reverse the parity of the fibers to turn these bundles into
%Lie supergroups, we find in general that $\Pi E \ncong G_{0}\times\Pi V$.  On the other
%hand, it follows from our discussion in Section \ref{section-thom-group} below that these two Lie supergroups nevertheless possess isomorphic (super) representation rings.)
Here, the $\g_{0}$ directions correspond to differentiation along the base $G_{0}$, while the $\Pi V$ directions correspond to interior contraction along the fibers of $\Lambda^{*}(E^{*})$. In particular, for the coadjoint
representation $V = \g_{0}^{*}$, we have $E = T^{*}G_{0}$, and the Lie
superalgebra corresponding to $\Pi(T^{*}G_{0})$ is $\g_{0} \oplus \Pi\g_{0}^{*}$.  We will revisit this example in Sections \ref{Weyl-GKRS} and \ref{section-dirac-induction}.

\begin{remark}
Such bundles can be trivialized by left translation as $E \cong G_{0}\times V$.
However, when we reverse the parity of the fibers to turn these bundles into
Lie supergroups, we find in general that $\Pi E \ncong G_{0}\times\Pi V$.  On the other
hand, it follows from our discussion in Section \ref{section-thom-group} below that these two Lie supergroups nevertheless possess isomorphic super representation rings.
\end{remark}

When $G$ is a Lie supergroup, the terms ``$G$-module''
and ``$G$-super\-module'' refer to representations of the associated Lie superalgebra $\g = \g_{0}\oplus\g_{1}$, for which the restriction to the even part $\g_{0}$ exponentiates to a representation of the underlying Lie group $G_{0}$.  In particular, a $G$-module carries representations of both the Lie group $G_{0}$ and the Lie superalgebra $\g = \g_{0}\oplus \g_{1}$, such that
for $g\in G_{0}$ and $X\in \g$, their actions satisfy $g \circ X \circ g^{-1} = \Ad_{g} X$. In terms of associative superalgebras, a $G$-module is a representation
of the semi-direct tensor product $\dist(G_{0}) \semitensor U(\g)$, with multiplication
\begin{equation}\label{eq:semitensor}
(\phi\semitensor X)(\psi\semitensor Y) = (\phi\conv\psi) \semitensor X \Ad_{\phi}Y
\end{equation}
for $\phi,\psi\in \dist(G_{0})$ and $X,Y\in U(\g).$
To further reconcile the $\dist(G_{0})$-action with the $U(\g_{0})$ action, we note
that there is an algebra homomorphism $\varepsilon: U(\g_{0}) \to \dist(G_{0})$ given on the generators $X\in \g$ by
\begin{equation}\label{eq:epsilon}
\varepsilon: X \mapsto \frac{d}{dt}\Bigr|_{t=0}\delta_{\mathrm{exp}(tX)},
\end{equation}
identifying the universal enveloping algebra with derivations at the identity (see \cite{AM}). We can now define the Lie supergroup version of the convolution algebra as the tensor product over $U(\g_{0})$:
\begin{equation}\label{eq:supergroup-convolution}
\dist(G) := \bigl( \dist(G_{0})\semitensor U(\g) \bigr)\,/\,\bigl( \varepsilon(X)\semitensor 1 = 1\semitensor X \text{ for $X\in U(\g_{0})$} \bigr).
\end{equation}
This algebra does indeed correspond to distributions on the Lie supergroup. For instance, for the purely odd Lie supergroup with $G_{0} = 1$ and $\g = 0 \oplus \Pi V$, it gives the ring of functions $\dist(G) \cong \Lambda^{*}(V)$.
%**********
%which we use to define the super representation ring:

\begin{definition}
Let $G$ be a Lie supergroup with underlying even Lie group $G_{0}$ and associated Lie superalgebra $\g = \g_{0} \oplus \g_{1}$.  Its
super representation ring is the graded ring with components
\begin{align*}
%	R(G) := R\bigl( \dist(G) \bigr), \qquad 
%	%R\bigl( \,\dist(G) \semitensor U(\g) \, \bigr), \\
	SR^{-n}(G) := SR^{-n}\bigl( \dist(G) \bigr) =
	SR \bigl( \dist(G) \tensor \Cl(n) \bigr),
	%SR^{-n}\bigl( \, \dist(G) \semitensor U(\g) \, \bigr),
\end{align*}
constructed from supermodules carrying compatible $G_{0}$- and $\g$-actions. 
\end{definition}

\section{Twisted representation rings}
\label{twisted}

\subsection{Twistings of Lie groups}
\label{twist-group}

Let $G$ be a compact (connected) Lie group.
%For certain cohomology classes
%$\tau\in H^{2}(G;\Z)$, we can construct a twisted version
%of the representation ring $R^{\tau}(G)$, consisting of virtual, finite dimensional, complex projective representations of $G$ with projective cocycle $\tau$.  For the trivial twisting $\tau = 0$, we recover the usual representation ring $R^{0}(G) = R(G)$.
%More specifically,
Recall that a projective unitary representation of $G$ on a complex Hilbert space $V$ is a continuous group homomorphism $r: G \to \PU(V) = \U(V)/\U(1)$.
%Such  projective representations are classified by the pullback $r^{*}\,c_{1}\bigl(\U(V)\bigr)\in H^{2}(G;\Z)$ of the Chern class of $\U(V)$ viewed as a principal $\U(1)$-bundle over $\PU(V)$.
Viewing $\U(V)$ as a principal $\U(1)$-bundle over $\PU(V)$,
pulling back this bundle to $G$ gives us a central extension $\widetilde{G} = r^{*}\U(V)$ of $G$, and the projective representation $r$ lifts to an actual representation $\tilde{r}$ on $\widetilde{G}$. This gives the following commutative diagram:
$$
\begin{CD}
  1 @>>> S^{1}      @>>> \widetilde{G} @>>> G @>>> 1 \\
  @. @VV\cong V @VV\tilde{r}V @VVrV @.\\
  1 @>>> \U(1) @>>> \U(V) @>>> \PU(V) @>>> 1
\end{CD}
%\begin{CD}
%  S^{1}      @>\cong>>   \U(1) \\
%  @VVV            @VVV  \\
%  \widetilde{G}  @>\tilde{r}>> \U(V)  \\
%  @VVV            @VVV  \\
%  G          @>r>> \PU(V)
%\end{CD}
$$
The central extension $\widetilde{G}$ is called the \emph{cocycle} of the projective representation, and it is classified topologically by its Chern class
$$c_{1}(\widetilde{G}) = c_{1}\bigl( r^{*}\U(V) \bigr) = r^{*} c_{1}\bigl( \U(V)\bigr) \in H^{2}(G;\Z).$$
%which is once again the cocycle of the projective representation.
However, the Chern
class itself does not capture the group structure of this extension. To see
the group structure topologically, we must ``deloop'' the groups and work instead with the corresponding fibration of classifying spaces:
\begin{equation}
\begin{CD}\label{eq:BG}
  S^{1} @>>> ES^{1} @>>> BS^{1}\\
  @VVV            @VVV   @VVV\\
  \widetilde{G} @>>> E\widetilde{G} @>>> B\widetilde{G} \\
  @VVV            @VVV  @VVV \\
  G @>>> EG @>>> BG
\end{CD}
\end{equation}
%Considering the transgressions for the fibrations along the top, bottom, left and right sides of (\ref{eq:BG}), we obtain the following commutative diagram in cohomology:
%$$\begin{CD}
%H^{1}(S^{1}) @>d_{2}^{ES^{1}}>> H^{2}(BS^{1}) \\
%@VVd_2^{\widetilde{G}}V @VVd_{3}^{B\widetilde{G}}V \\
%H^{2}(G) @>d_{3}^{EG}>> H^{3}(BG)
%\end{CD}$$
%Recalling that $BS^{1}\simeq \C P^{\infty}$ is a $K(\Z,2)$, the bundle $B\widetilde{G}$ on the right side of (\ref{eq:BG}) is
%determined topologically by its characteristic class in $H^{3}(BG;\Z)$.
Since $S^{1}$ is abelian, we can choose a $BS^{1}$, up to homotopy, which is a topological group and for which the universal fibration $S^{1} \to ES^{1} \to BS^{1}$ on the top of (\ref{eq:BG}) is a central extension. A circle bundle is determined by a continuous map $f: G \to BS^{1}$, and in order to pull back the group structure on $ES^{1}$ to $\widetilde{G} = f^{*}ES^{1}$, we require that $f$ be a group homomorphism. In that case, ``delooping'' $f$ gives a map of classifying spaces $Bf: BG \to BBS^{1} = K(\Z,3)$. This map determines the bundle $B\widetilde{G}$ on the right side of (\ref{eq:BG}), which is completely characterized up to isomorphism by the class $[Bf] \in [BG,K(\Z,3)] = H^{3}(BG;\Z)$.

Considering the transgressions for the fibrations along the top, bottom, left and right sides of (\ref{eq:BG}), we obtain the following commutative diagram in cohomology:
$$\begin{CD}
H^{1}(S^{1}) @>d_{2}^{ES^{1}}>> H^{2}(BS^{1}) \\
@VVd_2^{\widetilde{G}}V @VVd_{3}^{B\widetilde{G}}V \\
H^{2}(G) @>d_{3}^{EG}>> H^{3}(BG)
\end{CD}$$
%Recalling that $BS^{1}\simeq \C P^{\infty}$ is a $K(\Z,2)$, the bundle $B\widetilde{G}$ on the right side of (\ref{eq:BG}) is
%determined topologically by its characteristic class in $H^{3}(BG;\Z)$.
Letting $\theta$ 
be the generator of $H^{1}(S^{1};\Z)$, and $u = d_{2}^{ES^{1}} \theta$ the generator of $H^{2}(BS^{1};\Z)$, the characteristic class $[Bf]$ of the bundle $B\widetilde{G}$ is the transgression $d_{3}^{B\widetilde{G}}u = d_{3}^{B\widetilde{G}}d_{2}^{ES^{1}}\theta$.
On the other hand, we can also write this class 
%in terms of the characteristic class of $\widetilde{G}$
as
$d_{3}^{EG}c_{1}(\widetilde{G})
= d_{3}^{EG}d_{2}^{\widetilde{G}}\theta$.
%, where $\tau = c_{1}(\widetilde{G}) = d_{2}\theta \in H^{2}(G;\Z)$.

%\footnote{To construct this isomorphism $\Ext(G,S^{1}) \to H^{3}(BG)$ explicitly, recall that a circle bundle is given by a continuous map $f : G \to BS^{1}$. Since $S^{1}$ is abelian, we can choose a classifying space $BS^{1}$, up to homotopy, which is a topological group and for which the universal fibration $S^{1} \to ES^{1} \to BS^{1}$ is a central extension. In order to pull back the group structure on $ES^{1}$ to $\widetilde{G} = f^{*}ES^{1}$, we require that $f$ be a group homomorphism, in which case it induces a map of classifying spaces $BG \to BBS^{1} = K(\Z;3)$, which determines a class in $H^3(BG;Z) = [BG,K(\Z,3)]$.}

\begin{proposition}\label{prop-classify}
A circle bundle $S^{1}\to\widetilde{G}\to G$ admits the structure of a group extension precisely when the characteristic class $c_{1}(\widetilde{G})$ is in the kernel of the differential
$$d_{2} : H^{2}(G;\Z)  \longrightarrow E_{2}^{1,2}(EG;\Z) = H^{2}\bigl(BG;\,H^{1}(G;\Z)\bigr)$$ in the spectral sequence for the universal fibration $G\to EG\to BG$. In addition, the map
$$d_{3}: \mathrm{Ker}\,d_{2} \hookrightarrow H^{2}(G;\Z) \to H^{3}(BG;\Z)$$
is an isomorphism, and isomorphism classes of group extensions $S^{1}\to\widetilde{G}\to G$ correspond to elements of $H^{3}(BG;\Z)$
via the transgressions $d_{3}\,c_{1}(\widetilde{G})\in H^{3}(BG;\Z)$ of their
characteristic classes.
\end{proposition}

\begin{proof}
In \cite{AS1}, Atiyah and Segal discuss such extensions and prove that $\Ext(G,S^{1}) \cong H^{3}(BG;\Z)$. In summary, they show that $\Ext(G,S^{1}) \cong H^{2}(BG;S^{1})$ (see \cite{M1,M2}), and then they use the coefficient sequence $\Z \to \R \to S^{1}$ to establish the isomorphism $H^{2}(BG;S^{1}) \cong H^{3}(BG;\Z)$
for compact $G$.

Let us now examine the Serre spectral sequence for the universal fibration $G\to EG\to BG$, where $EG$ is contractible and thus its reduced cohomology $\widetilde{H}^{*}(EG)$ is trivial.
We first note that $H^{1}(BG)$ must vanish, which gives us $E_{2}^{1,1} = 0$.
It follows that $H^{3}(BG)$ survives past the $E_{2}$ stage,
and to kill it at the $E_{3}$ stage, the map $d_{3}$ must be an isomorphism from $\mathrm{Ker}\,d_{2}\subset H^{2}(G)$ to $H^{3}(BG)$.
We therefore obtain the exact sequence
\begin{equation}\label{eq:exact-sequence}
\begin{CD}
0 @>>> H^{3}(BG) @>d_{3}^{-1}>> H^{2}(G) @>d_{2}>> E_{2}^{1,2}
@>d_{2}>> H^{4}(BG),
\end{CD}\end{equation}
where the inverse transgression is the looping map
$d_{3}^{-1} : [BG,K(\Z,3)] \to [G, K(\Z,2)]$, induced by the forgetful map taking a central extension to its underlying topological circle bundle.
The Chern class of a central extension $\widetilde{G}$ is then the image
$c_{1}(\widetilde{G}) = d_{3}^{-1} ( d^{ B\widetilde{G}}_{3} u )$, and thus $d_{2}\,c_{1}(\widetilde{G})$ vanishes.
%In general, suppose we have a fibration $F \to E \to X$, where the cohomology $H^{*}(F)$ of the fiber is freely generated.  We can then classify this bundle by looking at the transgressions of those generators, giving us characteristic classes in $H^{*}(X)$.  Now consider the looped fibration $\Omega F \to \Omega E \to \Omega X$. The cohomology $H^{*}(\Omega F)$ of the fiber must once again be freely generated, and the generators of $H^{*}(F)$ are the transgressions of the generators of $H^{*}(\Omega F)$ in the spectral sequence for the path fibration $\Omega F \to PF \to F$. If $\theta\in H^{*}(\Omega F)$ is a generator, and $u\in H^{*}(F)$ is its transgression, then the characteristic classes $c_{\theta}\in H^{*}(\Omega X)$ and $c_{u}\in H^{*}(X)$ are related by
%$D c_{\theta} = c_{u}$, where $D$ is the differential in the spectral sequence for the path fibration $\Omega X \to PX \to X$.
\end{proof}

\begin{example}
For the torus $T^{2}=S^{1}\times S^{1}$, we have
$H^{2}(T^{2};\Z) = \Z$, and so it admits nontrivial circle bundles
$S^{1}\to \widetilde{T}^{2}\to T^{2}$. However, for the classifying space we
have
$H^{3}(BT^{2};\Z) = 0$, and thus $T^{2}$ does not admit any nontrivial
$S^{1}$ group extensions.  The bundle $\widetilde{T}^{2}$ therefore does
not admit a group structure. In this case, the exact sequence (\ref{eq:exact-sequence}) becomes
\begin{equation*}
0 \xrightarrow{\phantom{d_{3}^{-1}}} 0 \xrightarrow{d_{3}^{-1}} \Z \xrightarrow{\,\,\,d_{2}\,\,\,} \Z^{\oplus 4} \xrightarrow{\,\,\,d_{2}\,\,\,}  \Z^{\oplus 3}
\end{equation*}
and since $d_{2}$ is injective, the obstruction $d_{2}\,c_{1}(\widetilde{T}^{2})$
is nonzero.
\end{example}

\begin{definition}
Given a Lie group $G$, a \emph{twisting} $\tau$ is a central extension $$\tau = \{ 1 \to S^{1} \to \widetilde{G} \to G \to 1 \}$$
The $\tau$-twisted $G$-equivariant $K$-theory of a point, $\K{\tau}{}{G}$, is the Grothendieck group of isomorphism classes of finite dimensional projective representations of $G$ with cocycle $\tau$.
\end{definition}

This twisted $K$-group depends, up to noncanonical isomorphism, only on the isomorphism class of the twisting,  the class
$[\tau] = c_{1}(\widetilde{G})\in H^{2}(G;\Z)$, which by Proposition \ref{prop-classify} must satisfy $d_{2}[\tau] = 0$. Given two cocycles $\tau_{1}$, $\tau_{2}$ corresponding to central extensions $\widetilde{G}_{1}, \widetilde{G}_{2}$ respectively, their sum $\tau_{1}+\tau_{2}$ corresponds to the tensor product extension $\widetilde{G}_{1} \otimes \widetilde{G}_{2}$. Taking the tensor product of two projective representations therefore adds their cocycles. In particular, the twisted $K$-groups are not rings.
%The $\widetilde{G}$-equivariant $K$-theory of a point decomposes as
%$$\K{}{}{\widetilde{G}} = {\bigoplus}_{k\in\Z}\K{}{}{\widetilde{G}}_{k},$$ according to the action of the central $S^{1}$.
The twisted equivariant $K$-theory can then be constructed as the component $\K{\tau}{}{G} \cong \K{}{}{\widetilde{G}}_{1}$ of virtual $\widetilde{G}$-modules on which the central $S^{1}$ acts by complex multiplication.
%Furthermore, we have $\K{k\tau}{}{G} = \K{}{}{\widetilde{G}}_{k}$. The zero component consists of virtual $\widetilde{G}$-modules on which the central $S^{1}$ acts trivially, i.e., lifts of virtual $G$-modules, giving us $\K{}{}{G} \cong \K{}{}{\widetilde{G}}_{0}$. 
%The twisted groups $\K{\tau}{}{G}$ are not themselves rings; taking the tensor product of two projective representations adds their cocycles.
%However, taken in their entirety, the direct sum of all the twisted representation groups $\K{*}{}{G} = \bigoplus_{\tau}\!\K{\tau}{}{G}$ is an $H^{3}(BG;\Z)$-graded ring over $R(G)$, corresponding when $G$ is semisimple to the equivariant $K$-theory $\K{}{}{\bar{G}}$ for the universal cover $\bar{G}$ of $G$.

\begin{example}
Let $G = \SO(n)$ for $n\geq 3$. There are two
types of representations of $\SO(n)$ corresponding to the two elements of $H^{2}(\SO(n);\Z) \cong \Z_{2}$. The bosonic or integer spin representations are true $\SO(n)$-modules and comprise $\K{0}{}{\SO(n)}$.  We also have projective representations of $\SO(n)$.  In the simplest case, the Lie group isomorphism $\PU(2) = \U(2) / \U(1) \cong \SU(2) / \{\pm 1\} \cong \SO(3)$
yields a nontrivial projective representation of $\SO(3)$ on $\C^{2}$, whose
cocycle gives the nontrivial element $[\tau]\in H^{2}(\SO(3);\Z)$. In general, the fermionic or half-integer spin representations are projective representations of $\SO(n)$
whose cocycle gives the nontrivial element
$[\tau]\in H^{2}(\SO(n);\Z)$, and which comprise $\K{\tau}{}{\SO(n)}$.  These fermionic representations
can be constructed as true representations of the Lie algebra $\so(n)$. To realize them as representations of a Lie group,
we must lift to the universal cover $\Z_{2}\to\Spin(n)\to\SO(n)$. Then $\Spin(n)$-modules can be characterized as either bosonic or fermionic depending on whether the nontrivial central element $-1\in\Z_{2}\subset\Spin(n)$ acts by $+1$ or $-1$ respectively.
To phrase this construction in the terms of the preceding paragraph, we must modify this argument slightly to instead consider even and odd representations of the extension $S^{1}\to\Spin^{c}(n)\to\SO(n)$, which can be constructed as
$\Spin^{c}(n) = \Spin(n) \times_{\Z_{2}} S^{1}$.
\end{example}

For twistings by torsion cocycles, the projective representations of the Lie group are nevertheless true representations of the Lie algebra. This is because the projective representations lift to true representations
of a finite covering of the Lie group. Such finite covers
are locally indistinguishable from the original Lie group,
and so are not detected by the Lie algebra.

\subsection{Twistings of Lie algebras}
\label{twist-algebra}

We can also consider projective representations at the Lie algebra level. A projective unitary representation of a Lie algebra $\g$ on a complex
Hilbert space $V$ is a Lie algebra homomorphism $r: \g \to \mathfrak{pu}_{V} = \mathfrak{u}_{V} / \mathfrak{u}_{1}$, where $\mathfrak{u}_{1} = i\,\R\,\mathrm{Id}$. We can construct a splitting of the extension $\mathfrak{u}_{1}\to\mathfrak{u}_{V}\to\mathfrak{pu}_{V}$ by choosing a complement of $\mathfrak{u}_{1}$ in $\mathfrak{u}_{V}$ and identifying it with $\mathfrak{pu}_{V}$, giving us an isomorphism $\mathfrak{u}_{V} \cong \mathfrak{pu}_{V} \oplus \mathfrak{u}_{1}$.  With respect to this splitting, the bracket is given by
\begin{equation}\label{eq:u1-extension}
	[\, A \oplus a,\,B \oplus b\, ]_{\mathfrak{u}} =
	[A,B]_{\mathfrak{pu}} \oplus i\,\omega(A,B),
\end{equation}
for $A,B \in \mathfrak{pu}_{V}$ and $a,b\in \mathfrak{u}_{1}$.  Here,
the coefficient of the $\mathfrak{u}_{1}$-term is given by an anti-symmetric bilinear form $\omega : \mathfrak{u}_{V} \otimes \mathfrak{u}_{V} \to \R$, or in other words an element $\omega\in\Lambda^{2}(\mathfrak{pu}_{V}^{*})$.
As a consequence of the Jacobi identity, the 2-form $\omega$ satisfies the cocycle condition
\begin{equation}\label{eq:cocycle}
	\omega(A,[B,C]) = \omega([A,B],C) + \omega(B,[A,C]),
\end{equation}
%is invariant under the coadjoint action of $\mathfrak{u}_{V}$, satisfying
%\begin{equation}\label{eq:ad-invariance}
%\bigl(\ad^{*}_{C}\omega\bigr) (A,B) = \omega\bigl( [C,A], B %\bigr) + \omega \bigl( A, [C,B] \bigr) = 0.
%\end{equation}
or in other words, the form $\omega$ is closed, satisfying $d\omega = 0$ with respect to the Lie algebra cohomology differential $d : \Lambda^{*}(\mathfrak{pu}_{V}^{*})\to
\Lambda^{*+1}(\mathfrak{pu}_{V}^{*})$.  Furthermore, if we choose another splitting $\mathfrak{u}_{V} \cong \mathfrak{pu}_{V} \oplus \mathfrak{u}_{1}$, we obtain another cocycle which differs from $\omega$ by an exact form. Thus, $\omega$ determines a class in the Lie algebra cohomology $H^{2}(\mathfrak{pu}_{V})$ which is independent of the choice of splitting.
The \emph{cocycle} of a projective Lie algebra representation is the pullback $\tau = r^{*}(\omega)$, determining a class $[\tau] \in H^{2}(\g)$.

\begin{remark}If $V$ is finite dimensional, then we have a canonical inclusion
$\mathfrak{pu}_{V}\hookrightarrow\mathfrak{u}_{V}$ as the traceless elements.
This gives us a canonical splitting $\mathfrak{u}_{V} \cong \mathfrak{pu}_{V} \oplus \mathfrak{u}_{1}$ with respect to which the cocycle $\omega$ vanishes.  Thus, every projective unitary Lie algebra representation on a finite dimensional vector space lifts trivially to an actual representation. We can also see this from the point of view of Lie algebra cohomology, as $H^{2}(\mathfrak{pu}_{n})$ vanishes.  Indeed, the Lie algebra cohomology
$H^{2}(\g)$ vanishes for any finite dimensional semisimple Lie algebra $\g$.
\end{remark}

Given a projective Lie algebra representation $r:\g\to\mathfrak{pu}_{V}$, 
we can lift it to a vector space homomorphism $\tilde{r} : \g \to \mathfrak{u}_{V}$,
and in light of (\ref{eq:u1-extension}), the cocycle $\tau$ is a quantitative
measure of the failure of $\tilde{r}$ to be a Lie algebra homomorphism:
\begin{equation*}
	\tilde{r}(X)\,\tilde{r}(Y) - \tilde{r}(Y)\,\tilde{r}(X)
	= \tilde{r}\bigl([X,Y]\bigr) + i\,\tau(X,Y)\,\mathrm{Id}.
\end{equation*}
As in the Lie group case, we can pull back the extension, giving a commutative diagram
$$
\begin{CD}
	0 @>>> \R I @>>> \tilde{\g} @>>> \g @>>> 0\\
	@. @VV\cong V @VV\tilde{r}V @VVrV @. \\
	0 @>>> \mathfrak{u}_{1} @>>> \mathfrak{u}_{V} @>>> \mathfrak{pu}_{V} @>>> 0
\end{CD}
$$
where the central extension $\tilde{\g} = \g \oplus \R I$ of $\g$ has Lie algebra bracket
\begin{equation}\label{eq:extension-bracket}
	[\, X \oplus xI,\,Y \oplus yI \, ]_{\tilde{\g}} =
	[X,Y]_{\mathfrak{g}} \oplus \tau(X,Y)\,I,
\end{equation}
and the lift $\tilde{r}$ takes the central generator $I$ to $i\,\mathrm{Id}$.
Note that the cocycle condition (\ref{eq:cocycle}) for $\tau$
ensures that the Jacobi identity on $\g$ extends to $\tilde{\g}$.

Unlike in the Lie group case, \emph{every} Lie algebra cohomology 2-cocycle $\tau
%\in H^{2}(\g)
$ gives rise to a central extension $\R I\to\tilde{\g}\to\g$
via the bracket (\ref{eq:extension-bracket}).  In fact, the central $\mathfrak{u}_{1}$-extensions of $\g$ are classified up to isomorphism by $H^{2}(\g)$.
If $\g$ is the Lie algebra of a compact Lie group $G$, then a Lie algebra extension of
$\g$ given by the cocycle $\tau$ must satisfy two conditions in order for it to exponentiate to give a group extension of $G$. First, the cocycle must be integral, corresponding to a class $$[\tau]\in H^{2}(G;\Z) \to H^{2}(G;\R)\cong H^{2}(\g).$$ Second, the cocycle must be integrable, satisfying
$$d_{2}[\tau] = 0 \in \Lambda^{1}(\g^{*}) \otimes \Sym^{1}(\g^{*}) \cong E^{1,2}_{2}(EG)$$
in the Weil algebra for $\g$, viewed as a model for the de Rham cohomology of $EG$.

\begin{example}
The two-dimensional abelian Lie algebra $\mathfrak{t}^{2}$ is the Lie algebra of the torus $T^{2} = S^{1}\times S^{1}$ which we considered in the previous section.  We have $H^{2}(\mathfrak{t}^{2}) \cong \mathbb{R}$, and thus $\mathfrak{t}^{2}$ does indeed admit a central extension.  Letting $e,f$ be a basis for $\mathfrak{t}^{2}$, this extension has brackets $[e,f] = I$ and $[I,e] = [I,f] = 0$. This Lie algebra is the Heisenberg algebra, which has the standard infinite dimensional representation on the space of polynomials $\C[t]$ (or its $L^{2}$-completion), with action given by
$e \mapsto t$, $f \mapsto -i\,\partial_{t}$, $I \mapsto i\,\mathrm{Id}$,
from quantum mechanics.  We note that this extension does \emph{not} exponentiate to
give an extension of the Lie group $T^{2}$.
\end{example}

\begin{definition}
Given a Lie algebra $\g$ and a Lie algebra cohomology 2-cocycle $\tau$ representing a class $[\tau]\in H^{2}(\g)$, the $\tau$-twisted $\g$-equivariant $K$-theory of a point, $\K{\tau}{}{\g}$, is the Grothendieck group of isomorphism classes of finite dimensional projective representations of $\g$ with cocycle $\tau$.  
\end{definition}

As in the Lie group case, the twisted $K$-group depends, up to noncanonical isomorphism, only on the class $[\tau]\in H^{2}(\g)$. Taking the tensor product of two projective representations of $\g$ adds their cocycles, in this case adding them as 2-forms.
%Taking all the twistings together, we see that $\K{*}{}{\g}$ is a $H^{2}(\g)$-graded ring over $\K{}{}{\g}$, which corresponds to $\K{}{}{\bar{\g}}$ for the universal central extension $\bar{\g}$ of $\g$.
To introduce degree shifts, we must construct an associative algebra version of these twistings.  Given a cocycle $\tau\in H^{2}(\g)$, let $\tilde{\g}_{\tau}$ denote the central extension of $\g$ given by (\ref{eq:extension-bracket}).  A $\tau$-twisted representation of $\g$ is then a representation of $\tilde{\g}_{\tau}$ where $I$ acts by $i\,\mathrm{Id}$.  We thus construct the $\tau$-twisted universal enveloping algebra of $\g$ as the quotient
\begin{equation}\label{eq:tuea}
U_{\tau}(\g) := U(\tilde{\g}_{\tau})\,/\,( I = i\,1).
\end{equation}
The $\tau$-twisted representations of $\g$ are clearly representations of $U_{\tau}(\g)$, and we define:

\begin{definition}
Let $\g$ be a Lie algebra.  Given a 2-cocycle $\tau$ representing a class $[\tau]\in H^{2}(\g)$, the $\tau$-twisted super representation group of $\g$ is the graded group with components
\begin{equation}\label{eq:twisted-lie-superalgebra}
%\RR{\tau}{}(\g) := R\bigl( U_{\tau}(\g) \bigr), \qquad
  \SR{\tau}{-n}(\g) := SR^{-n}\bigl( U_{\tau}(\g) \bigr) = SR \bigl( U_{\tau}(\g) \tensor \Cl(n) \bigr),
\end{equation}
constructed from projective representations of $\g$ with cocycle $\tau$.
\end{definition}

Revisiting twistings of Lie groups, consider a twisting $\tau = \{ 1 \to S^{1} \to \widetilde{G} \to G \to 1 \}$. The corresponding extension $0 \to \R I \to \tilde{\g} \to \g \to 0$ of Lie algebras determines a Lie algebra cohomology 2-cocycle, which we also denote by $\tau$.
%class $\tau\in H^{2}(G;\Z)$ with vanishing transgression $d_{2}\tau = 0$.
%Such a class $\tau$ gives a cocycle in the Lie algebra cohomology $H^{2}(\g)$, which we also denote by $\tau$.  If $\widetilde{G}_{\tau}$ is
%the extension of $G$ with $c_{1}(\widetilde{G}) = \tau$, then its Lie algebra
%is the extension $\tilde{\g}_{\tau}$.
A projective representation of $G$ with cocycle $\tau$ is thus a projective representation of $\g$ with cocycle $\tau$, so we can describe a $\tau$-twisted representation of $G$ as a representation of $\widetilde{G}_{\tau}$, on which the corresponding $\tilde{\g}_{\tau}$-action satisfies $I = i\,\mathrm{Id}$. Recalling the algebra homomorphism $\varepsilon : U(\tilde{\g}_{\tau}) \to \dist(\widetilde{G}_{\tau})$ given by (\ref{eq:epsilon}), we define the twisted convolution algebra,
\begin{equation}\label{eq:twisted-cg}
\dist_{\tau}(G) := \dist(\widetilde{G}_{\tau}) \, / \, ( \varepsilon(I) = i\,1 ),
%\bigl( \dist(\widetilde{G}_{\tau}) \semitensor U(\tilde{\g}_{\tau}) \bigr) \, / \, \bigl(1 \semitensor I = i\,1\semitensor 1 \bigr)
%\cong
%\dist(\widetilde{G}_{\tau}) \semitensor U_{\tau}(\g),
\end{equation}
%where the semi-direct tensor product has multiplication given by (\ref{eq:semitensor}) for elements $\phi,\psi\in \dist(\widetilde{G}_{\tau})$ and $X,Y\in U(\tilde{\g}_{\tau})$.
which we use in the following definition:

\begin{definition}
Let $G$ be a compact Lie group.  Given a 2-cocycle $\tau$ representing a class $[\tau]\in H^{2}(G;\Z)$ with 
%vanishing transgression
$d_{2}[\tau] = 0$, the $\tau$-twisted super representation group of $G$ is the graded group with components
\begin{equation}\label{eq:twisted-lie-group}
%	\RR{\tau}{}(G) := R \bigl( \dist_{\tau}(G) \bigr), \qquad
	\SR{\tau}{-n}(G) := SR^{-n} \bigl( \dist_{\tau}(G) ) 
	= SR \bigl( \dist_{\tau}(G)\tensor \Cl(n) \bigr),
\end{equation}
constructed from projective representations of $G$ with cocycle $\tau$.
\end{definition}

\subsection{Twistings of Lie superalgebras and supergroups}
\label{twist-superalgebra}

The discussion of the previous section works equally well for Lie superalgebras $\g = \g_{0}\oplus \g_{1}$. A (even) extension $0 \to \R I \to \tilde{\g}_{\tau} \to \g \to 0$ is once again determined by a 2-cocycle $\tau$. In this case, the cocycle
is a supersymmetric 2-form in
\begin{equation}\label{eq:super-exterior}
\Lambda^{2}(\g_{0}^{*}\oplus \g_{1}^{*})
= \Lambda^{2}(\g_{0}^{*}) \oplus \Sym^{2}(\g_{1}^{*}) \oplus (\g_{0}^{*}\otimes \g_{1}^{*}),
\end{equation}
and it satisfies the supersymmetric version of the cocycle condition (\ref{eq:cocycle}):
$$
	\tau(A,[B,C]) = \tau([A,B],C) + (-1)^{|A|\,|B|}\,\tau(B,[A,C]),
$$
for homogeneous elements $A,B,C\in \g$. Such a cocycle determines an element $[\tau]\in H^{2}(\g)$ in the Lie superalgebra cohomology. The definitions (\ref{eq:tuea}) of the twisted univeral enveloping algebra and (\ref{eq:twisted-lie-superalgebra}) of the twisted super representation ring then carry over directly to Lie superalgebras.

\begin{example} Let $\g = 0 \oplus V$ be a purely odd Lie superalgebra.  Since the product of two odd elements is even, this $\g$ is supercommutative. 
%A projective cocycle in $H^{2}(\g)$ is then a $\g$-invariant form in $\Lambda^{2}(\g)$.
Consequently, every form is closed,
and recalling (\ref{eq:super-exterior}), a projective cocycle is a symmetric bilinear form
$b\in \Sym^{2}(V)$.  The extension $\tilde{\g}_{b}$ now satisfies
$$[v,w]_{\tilde{\g}} = [v,w]_{\g} \oplus b(v,w)\,I.$$
for $v,w\in V$, and the $b$-twisted universal enveloping algebra is then
$$U_{b}(\g) = T^{*}(V) / \bigl( v \cdot w + w \cdot v = i\,b(v,w) \bigr).$$
Up to replacing the coefficient $i$ on the right side with $-2$ (which we can do over $\C$ by scaling $V$ by $\sqrt{2i}$), this is the
definition (\ref{eq:clifford}) of the Clifford algebra $\ClR(V,b)$, which gives us an isomorphism
$$%\RR{b}{}(0\oplus V) \cong \RR{}{}\bigl(\ClR(V,b)\bigr), \qquad
  \SR{b}{-n}(0\oplus V) \cong \SR{}{-n}\bigl(\ClR(V,b)\bigr).$$
%of the super representation ring.
If $b$ is a \emph{non-degenerate} symmetric bilinear form on $V$, then $\SR{b}{*}(0\oplus V) \cong \SR{}{*+\dim V}(0).$
\end{example}

%\begin{example}
%Let $\g_{0}$ be a conventional even Lie algebra, and let $r:\g_{0}\to\End(V)$
%be a representation. We can then construct the Lie superalgebra $\g = \g_{0}\oplus \Pi V$ in which the bracket of any any two odd elements vanishes.  A symmetric bilinear form $b : V \tensor V \to \R$ on $V$ is $\g_{0}$-invariant if it satisfies
%$b(Xv,w) + b(v,Xw) = 0$.  Since the odd component $V$ acts trivially on itself here, a $\g_{0}$-invariant form is automatically $\g$-invariant, and thus gives an element
%$$b\in \bigl( \Sym^{2}(V^{*}) \bigr)^{\g_{0}} = \bigl( \Sym^{2}(V^{*}) \bigr)^{\g}
%\subset \bigl( \Lambda^{2}(\g^{*}) \bigr)^{\g} \cong H^{2}(\g).$$
%The form $b$ then defines a twisting, and restricting to the odd part, we have
%$U_{b}(V) \cong \ClR(V,b).$  For the entire Lie superalgebra, we then have
%\end{example}

%\subsection{Twistings of Lie supergroups}

Let $G$ be a Lie supergroup, with underlying even Lie group $G_{0}$ and Lie superalgebra $\g = \g_{0} \oplus \g_{1}$.  A twisting for $G$ is a pair $\tau = (\tau_{G_{0}},\tau_{\g})$, where $\tau_{G_{0}}$ is
a cocycle for the Lie group $G_{0}$ and $\tau_{\g}$ is a cocycle for the Lie superalgebra $\g$.  Letting $\tau_{\g_{0}}$ be the cocycle for $\g_{0}$ induced by $\tau_{G_{0}}$, we require the compatibility condition $\tau_{\g_{0}} = (\tau_{\g}) |_{\g_{0}}$.
%the compatibility condition that
%$d_{3}^{-1}[\tau_{G_{0}}]\mapsto [\tau_{\g}] $
%with respect to the maps:
%$$
%\begin{CD}
%	H^{3}(BG_{0};\Z) @. @. H^{2}(\g) \\
%	@A{d_{3}}A{\cong}A @. @VV{i^{*}}V \\
%	\mathrm{Ker}\,d_{2}  @>>>  H^{2}(G_{0};\Z)  @>>>  H^{2}(\g_{0})
%\end{CD}
%\begin{diagram}
%	H^{3}(BG_{0};\Z) &&& & H^{2}(\g) \\
%	\uTo>{\cong}<{d_{3}} &&& & \dTo>{i^{*}} \\
%	\mathrm{Ker}\,d_{2}&\rInto &H^{2}(G_{0};\Z) & \rTo & H^{2}(\g_{0})
%\end{diagram}
%$$
%recalling the transgressions $d_{2}$ and $d_{3}$ described in Section \ref{twist-group}.
Such a twisting determines extensions $\widetilde{G}_{0}$ of the Lie group $G_{0}$
and $\tilde{\g}$ of the Lie superalgebra $\g_{0}$, which
are compatible in that the Lie algebra associated to $\widetilde{G}_{0}$ is indeed the even component of $\tilde{\g}$.  These two extensions then combine to give a Lie supergroup central extension $\widetilde{G}_{\tau}$ of $G$.
%A projective representation of $G$ with cocycle $\tau$ is
%then a $\widetilde{G}_{\tau}$-supermodule on which the central generator $I$ acts by $i\,\mathrm{Id}$.  We recall from Section \ref{supergroups}
%that a $\widetilde{G}_{\tau}$-supermodule carries compatible actions
%of both the associated Lie superalgebra $\tilde{\g}_{\tau}$ and the underlying even Lie group $\widetilde{(G_{0})}_{\tau_{G_{0}}}$.
%%, which intertwine according to the identity $g \circ X \circ g^{-1} = \Ad_{g} X$ for $g \in \widetilde{(G_{0})}_{\tau_{G_{0}}}$ and $X\in\g$.
The definitions (\ref{eq:twisted-cg}) of the twisted convolution algebra and (\ref{eq:twisted-lie-group}) of the twisted super representation ring then carry over directly to the Lie supergroup case.

\section{The Thom isomorphism}
\label{thom}

In this section we consider an algebraic Thom isomorphism between the twisted super representation rings of a Lie algebra or compact Lie group and a related Lie superalgebra or Lie supergroup.  This serves as an instructive example involving twistings for both an even Lie group and the odd part of a Lie superalgebra. It is also an essential result which we will use in Sections \ref{Weyl-GKRS} and \ref{section-dirac-induction} below.

\subsection{Lie algebra version}
\label{thom-algebra}

%First we consider the Lie algebra/superalgebra version of the Thom isomorphism.
Let $\g$ be a finite dimensional Lie algebra, and let $r : \g\to\End(V)$ be a real finite dimensional representation.  We can then construct the Lie superalgebra $\g\oplus \Pi V$, with $\g$ for the even component and $V$ for the odd component, where the bracket of any two odd elements vanishes. (We write $\Pi V$ here as a reminder that $V$ is to be treated as an odd vector space.)
The universal enveloping algebra of the Lie superalgebra $\g\oplus \Pi V$ is the semi-direct tensor product
$$U(\g\oplus \Pi V) \cong U(\g) \semitensor \Lambda^{*}(V)$$
with multiplication given by
\begin{equation}\label{eq:semitensor1}
(X \semitensor v) (Y \otimes w) = XY \semitensor v\wedge w +
Y \semitensor v\wedge(r_{X}w)
\end{equation}
for $X,Y\in U(\g)$ and $v,w \in\Lambda^{*}(V)$, where we have extended the representation $r$ on $V$ to the exterior algebra $\Lambda^{*}(V)$ as a derivation.  Both $U(\g)$ and $U(\Pi V) \cong \Lambda^{*}(V)$ inject as algebras into the semi-direct tensor product as $U(\g) \semitensor 1$ and $1\semitensor\Lambda^{*}(V)$, and the multiplication (\ref{eq:semitensor1}) is defined so that $$(X\semitensor 1)(1\semitensor v) - (1\semitensor v)(X\semitensor 1)
= 1\semitensor r_{X}v$$
for $X\in \g$ and $v\in V$.

We say that a symmetric bilinear form $b$ on $V$ is
$\g$-invariant if it satisfies the identity $$(r_{X}b)(v,w) = b(r_{X}v,w) + b(v,r_{X}w) = 0$$ for all $X\in\g$ and $v,w\in V$. Since the odd component $\Pi V$ acts trivially on itself in the Lie superalgebra $\g\oplus \Pi V$, a $\g$-invariant form on $V$ is automatically $\g\oplus \Pi V$-invariant, and since invariant forms are closed, the form $b$ is a cocycle
% and we have
%$$b\in \bigl( \Sym^{2}(V^{*}) \bigr)^{\g} = \bigl( \Sym^{2}(V^{*}) \bigr)^{\g\oplus\Pi V}
%\subset \bigl( \Lambda^{2}(\g\oplus\Pi V)^{*} \bigr)^{\g\oplus\Pi V} \cong H^{2}(\g\oplus\Pi V).$$
defining a twisting for $\g\oplus \Pi V$. Restricting to the odd component we recall from our example at the end of Section \ref{twist-superalgebra} that $U_{b}(\Pi V) \cong \ClR(V,b)$.  This twisting affects only the odd component of the Lie superalgebra, and thus the twisted universal enveloping algebra for the full Lie superalgebra is
$$U_{b}(\g \oplus\Pi V) \cong U(\g) \semitensor \ClR(V,b),$$
where the multiplication on the semi-direct tensor product is once again given by (\ref{eq:semitensor1}), but this time taking $v,w\in \ClR(V,b)$, replacing wedge products with Clifford products, and extending $r$ to $\ClR(V,b)$ as a derivation with respect to the Clifford product.

Before stating the Thom isomorphism, we first recall the standard Lie algebra isomorphism $\so(V) \cong \spin(V)\subset\ClR(V)$, as described for example in \cite{LM}.

%Given a vector space $V$ with a symmetric bilinear form $b:V\otimes V \to \R$,
%we compare the Clifford algebra $\ClR(V)$ with the exterior algebra $\Lambda^{*}(V)$.
%We define a Clifford action $c : \ClR(V) \to \End(\Lambda^{*}(V))$ by $c(v) = e(v) + i(v^{*})$, where $e(v)$ is exterior multiplication, and $i(v^{*})$ is interior contraction by the element $v^{*}\in V^{*}$ satisfying $v^{*}(w) = b(v,w)$ for all $v,w\in V$. The Chevalley map $\ch : \ClR(V) \to \Lambda^{*}(V)$ is then given by $\ch(\omega) = c(\omega) 1$, which is an isomorphism of left $\ClR(V)$-modules. 

\begin{lemma}
	Let $V$ be a finite dimensional inner product space with orthonormal basis $\{e_{1},\ldots,e_{n}\}$.
	The extension of any $A\in\so(V)$ to the Clifford algebra $\ClR(\g)$
	as a derivation satisfying
	$$A (\xi\eta) = (A \xi)\,\eta + \xi\,(A\eta)$$
	for $\xi,\eta\in\ClR(V)$ can be quantized as the inner derivation given
	by bracketing with
	$$\widetilde{A} = -\frac{1}{4}\sum_{i=1}^{n} e_{i}\cdot(A e_{i}) \in \spin(V) \subset \ClR(V).$$
	Furthermore, the operator $A \mapsto \widetilde{A}$ is a Lie algebra homomorphism: $[\widetilde{A,B}] = \widetilde{A}\cdot\widetilde{B} - \widetilde{B}\cdot\widetilde{A}$.
\end{lemma}

\begin{proof}
	We need only verify that the operators $A$ and $[\widetilde{A},\,\cdot\,]$ agree
	on the generators
	$\xi\in V\subset \ClR(V)$, as the rest follows from the derivation property.
	Let $\xi = \sum_{i=1}^{n}\xi^{i}\,e_{i}$.  Then we have
	\begin{equation*}\begin{split}
		\widetilde{A}\cdot\xi - \xi\cdot\widetilde{A}
		&= -\frac{1}{4}\sum_{i=1}^{n} \bigl(
		-2b(e_{i},\xi) Ae_{i} -(-2)
		e_{i}\,b(Ae_{i},\xi) 
		\bigr) \\
		&= \frac{1}{2}\sum_{i=1}^{n} \bigl(
		b(e_{i},\xi)\,Ae_{i} +
		b(e_{i},A\xi)\,e_{i}
		\bigr)
		= \frac{1}{2}\sum_{i=1}^{n} \bigl(
		\xi^{i}\,Ae_{i} +
		(A \xi)^{i}\,e_{i} 
		\bigr) = A\xi,
	\end{split}\end{equation*}
	where $A\xi = \sum_{i=1}^{n}(A\xi)^{i}e_{i}$. For the Lie algebra property,
	we have
	$$\bigl[[\widetilde{A},\widetilde{B}],\xi\bigr] = \bigl[\widetilde{A},[\widetilde{B},\xi]\bigr] - \bigl[\tilde{B},[\widetilde{A},\xi]\bigr] = AB\xi - BA\xi = [A,B]\,\xi$$
	by the Jacobi identity. Since the center of $\ClR(V)$ consists of only 	
	scalar multiples of the identity, we have $[\widetilde{A,B}] = \widetilde{A}\cdot\widetilde{B} - \widetilde{B}\cdot\widetilde{A}$ up to addition by a scalar. In this case, the scalar
	terms vanish, as
	$$\widetilde{A}\cdot\widetilde{B} - \widetilde{B}\cdot\widetilde{A}
	= -\frac{1}{4} \sum_{i=1}^{n}\bigl( Ae_{i}\cdot Be_{i} + e_{i}\cdot AB e_{i} \bigr)
	= -\frac{1}{4} \sum_{i=1}^{n} e_{i}\cdot(AB -BA)e_{i} 
	= [\widetilde{A,B}],$$
	using the fact that $A$ is in $\so(V)$ (and that $V$ is finite dimensional) to obtain the second equality.
\end{proof}

\begin{proposition}[Thom Isomorphism]
If $b$ is a \emph{non-degenerate} $\g$-invariant symmetric bilinear form on $V$, then there is an isomorphism of unital associative algebras
$$U(\g) \semitensor \ClR(V,b) \xrightarrow{\cong} U(\g) \otimes \ClR(\dim V)$$
which induces an additive group isomorphism
\begin{align*}
%R(\g) \xrightarrow{\cong} \RR{b}{}( \g \oplus \Pi V ), \qquad
SR^{*}(\g) \xrightarrow{\cong} \SR{b}{*+\dim V}( \g \oplus \Pi V )
\end{align*}
between the super representation ring of $\g$ and the $b$-twisted super representation group of $\g \oplus \Pi V$.
\end{proposition}

\begin{proof}
Choosing an orthonormal basis $\{e_{1},\ldots,e_{\dim V}\}$ for $V$, we can identify $\ClR(V)$ with $\ClR(\dim V)$.  So, our goal is to construct an algebra isomorphism
$$f : U(\g) \semitensor \ClR(V) \to U(\g) \otimes \ClR(V).$$
Since the inner product $b$ is $\g$-invariant, the representation $r$ is orthogonal, giving a Lie algebra homomorphism $r : \g \to \so(V)$.
Using the above proposition to quantize this action, we obtain a Lie algebra homomorphism $\tilde{r} : \g \to \spin(V) \subset \ClR(V)$ (see \cite{KS}, \cite{L1}), which in turn lifts to the universal enveloping algebra to give a homomorphism $\tilde{r} : U(\g) \to \ClR(V)$ of associative (super) algebras.  Combining these two algebras, we get a tensor product representation $s : \g \to U(\g) \otimes \ClR(V)$,
\begin{equation}\label{eq:r}
	s(X) = X \otimes 1 + 1 \otimes \tilde{r}(X)
\end{equation}
for $X\in\g$, and we clearly have $s(X)\,s(Y) - s(Y)\,s(X) = s([X,Y])$.
This representation then extends to an algebra homomorphism $s : U(\g) \to U(\g) \otimes \ClR(V)$.  We can now construct $f$ as
$$f : X \semitensor v \mapsto s(X) (1 \otimes v) = X \otimes v + 1 \otimes \tilde{r}(X)\cdot v$$
for $X\in \g$ and $v\in\ClR(V)$.  

We observe that both $U(\g)$ and $\ClR(V)$ inject as algebras into their semi-direct tensor product $U(\g) \semitensor \ClR(V)$ as 
$U(\g) \semitensor 1$ and $1 \semitensor \ClR(V)$ respectively. Furthermore, their mixed products satisfy
\begin{equation}\label{eq:identity}
(X\semitensor 1)(1\semitensor v) - (1\semitensor v)(X\semitensor 1) = 1\semitensor r_{X}v
\end{equation}
for $X\in \g$ and $v\in V$.  To verify that $f$ is an algebra homomorphism, we check that it respects (\ref{eq:identity}): \begin{equation*}\begin{split}
	f(X\semitensor 1)&\, f(1\semitensor v) - f(1\semitensor v)\,f(X\semitensor 1)\\
	&=
	X \otimes v
	+ 1 \otimes \tilde{r}(X) \cdot v
	- X \otimes v
	- 1 \otimes v \cdot \tilde{r}(X)
	= f ( 1\semitensor r_{X} v ).
\end{split}\end{equation*}
Furthermore, this map is invertible, with inverse map
$$f^{-1} : U(\g) \otimes \ClR(V) \to U(\g) \semitensor \ClR(V)$$
given by
$$f^{-1} : X \otimes v \mapsto X \semitensor v - 1 \semitensor \tilde{r}(X) \cdot v.$$
Again, this map $f^{-1}$ is an algebra homomorphism, since we have
\begin{equation*}\begin{split}
	f^{-1}(X\otimes 1)&\, f^{-1}(1\otimes v)
	- f^{-1}(1\otimes v) \,f^{-1}(X\otimes 1)\\
	&=
	X \semitensor v - 1 \semitensor \tilde{r}(X) \cdot v
	+ 1 \semitensor r_{X}v
	- X \semitensor v
	+ 1 \semitensor v\cdot \tilde{r}(X)
	= 0,
\end{split}\end{equation*}
as we expect since $U(\g) \otimes 1$ commutes with $1 \otimes \ClR(V)$ in the standard tensor product. We see that $f \circ f^{-1} = f^{-1} \circ f = \mathrm{Id}$ acting on $X\in \g\subset U(\g)$ and $v\in\ClR(V)$, and thus these maps are inverses.

Pulling back $b$-twisted projective representations of $\g\oplus \Pi V$ to $\g$-Clifford supermodules via the algebra isomorphism $f$, we obtain an additive homomorphism of (twisted) super representation groups $f^{*} : \SR{}{-\dim V}(\g) \to \SR{b}{}(\g\oplus\Pi V)$ (see Section \ref{pullbacks} below for a discussion of pullbacks).  In fact, this map is an isomorphism, with inverse $(f^{*})^{-1} = (f^{-1})^{*}$.  Adding in degree shifts, we obtain our desired Thom isomorphism $\SR{}{*}(\g) \cong \SR{b}{*+\dim V}(\g\oplus\Pi V)$.
%We obtain the Thom isomorphism for the (twisted) representation groups by the same argument, but without the degree shifts.
\end{proof}

%We note that this twisted super representation ring can also be interpreted in terms of equivariant $K$-theory (with compact supports) as:
%$$\SR{b}{*}(\g \oplus \Pi V) \cong K_{\g}^{*}(V),$$
%in which case our Thom isomorphism becomes an isomorphism $K_{\g}^{*}(\mathrm{pt}) \cong K_{\g}^{*+\dim V}(V)$.

To construct this Thom isomorphism explictly, we recall that classes in $\SR{b}{\dim V}(\g\oplus\Pi V)$ correspond to $b$-twisted projective representations of $\g\oplus \Pi V$ with auxiliary supercommuting $\ClR(\dim V)$-actions, modulo those admitting odd involutions.  More precisely, we actually
want a Clifford action of $\ClR(0,\dim V)$ with generators $\{f_{1},\ldots,f_{\dim V}\}$
each squaring to $+1$, although in the complex case $\Cl(0,\dim V) \cong \Cl(\dim V)$.  Given
a $\g$-module $U$, we can then construct
$$f^{*} U = U \otimes \Cl(V),$$
where $\Cl(V) = \ClR(V) \otimes \C$ is the complex Clifford algebra.
The $b$-twisted $\g\oplus V$-action is given by
\begin{align}
	\label{eq:X}
	X(u\otimes \omega)
	&= X(u)\otimes \omega + 1 \otimes \tilde{r}(X) \cdot \omega, \\
	\label{eq:xi}
	v (u\otimes\omega) 
	&=  1 \otimes v \cdot \omega, \\
	\intertext{using the left action of $\ClR(V)$ on $\Cl(V)$,
	and the auxiliary $\ClR(0,\dim V)$-action uses the right action}
	\label{eq:fi}
	f_{i} (u\otimes\omega)
	&= u \otimes \grading(\omega) \cdot e_{i},
\end{align}
for $X\in \g$, $v \in V$, $u\in U$, $\omega\in \Cl(V)$.  Here,
the $f_{i}$ are Clifford generators squaring to $+1$, the $e_{i}$ are the corresponding orthonormal basis of $V$, and
$\grading : \Cl(V) \to \Cl(V)$ is the grading involution.
%Continuing the analogy to
%algebraic topology, the Thom class is $[\Cl(V)] \in \SR{b}{\dim V}(\g\oplus\Pi V)$, and the Euler class is its restriction to $\g$:
%$$e = \bigl[\Cl_{0}(V)\bigr] - \bigl[\Cl_{1}(V)\bigr] = \bigl[\Lambda^{\mathrm{even}}(V \otimes\C)\bigr]
%- \bigl[\Lambda^{\mathrm{odd}}(V \otimes \C)\bigr] \, \in \, SR(\g) \cong K_{\g}(\mathrm{pt}).$$
%(See Section \ref{pullbacks} below for a discussion of restrictions.)

\subsection{Lie group version}
\label{section-thom-group}

In the Lie group/supergroup case, we must consider an additional twisting.
Let $G$ be a compact Lie group, and let $V$ be a representation of $G$. Let $E$ be the vector bundle $E = (G \times G) \times_{\Delta G} V$ over $G$ with fibers $V$. Then $\Pi E$ is a Lie supergroup with underlying even Lie group $G$ and associated Lie superalgebra $\g \oplus \Pi V$ as described in Section \ref{supergroups}.

If $V$ is a finite dimensional real inner product space, we have the group extension
$$1 \to S^{1} \to \Spinc(V) \to \SO(V) \to 1$$
whose Chern class is the nontrivial element
$c_{1}(\Spinc(V)) \in H^{2}\bigl(\SO(V);\Z\bigr) \cong \Z_{2}.$
In our case, given a $G$-invariant inner product $b$ on $V$, then we have
a representation $r : G \to \SO(V)$, which we quantize by pulling back
the extension $\Spinc(V)$ via $r$:
$$
\begin{CD}
	1 @>>> S^{1} @>>> \widetilde{G} @>>> G @>>> 1 \\
	  @. @VV=V @VV\tilde{r}V @VVrV \\
	1 @>>> S^{1} @>>> \Spinc(V) @>>> \SO(V) @>>> 1
\end{CD}
%\begin{diagram}
%	1 & \rTo & S^{1} & \rTo & \widetilde{G} & \rTo & G & \rTo & 1 \\
%	  && \dTo>{=} && \dTo>{\tilde{r}} && \dTo>{r} &&\\
%	1 & \rTo & S^{1} & \rTo & \Spinc(V) & \rTo & \SO(V) & \rTo & 1
%\end{diagram}
$$
giving a continuous homomorphism $\tilde{r} : \widetilde{G} \to \Spinc(V) \subset\Cl(V)$.  Let the twisting $\tau$ be the central extension $\widetilde{G} = r^{*}\Spinc(V)$, which corresponds to a class
%This central $S^{1}$-extension $\widetilde{G}$ of $G$ is classified by its Chern class
\begin{equation}\label{eq:spinc-twisting}
[\tau] = c_{1}\bigl(\widetilde{G}\bigr)
= c_{1}\bigl( r^{*}\Spinc(V)  \bigr)
= r^{*} c_{1}\bigl(\Spinc(V)\bigr) \in H^{2}(G;\Z).
\end{equation}
%and thus $r$ is a projective representation of $G$ with cocycle $\tau$.
Exponentiating (\ref{eq:X}), the even $G$-component of the $b$-twisted $\Pi E$-action on $U \otimes \Cl(V)$ is given by
\begin{equation}\label{eq:g}
	g(u\otimes \omega) = g(u) \otimes \tilde{r}(g)\cdot\omega
\end{equation}
for $g\in \widetilde{G}$, $u\in U$, and $\omega\in \Cl(V)$. The odd
components act on $U \otimes \Cl(V)$ the same as they do in (\ref{eq:xi}) and (\ref{eq:fi}) above.  In order that the action (\ref{eq:g}) descend to $G$,
we see that $U$ must be a projective representation of $G$ with the opposite cocycle $-\tau$.  However, since $[\tau]$ is the pullback of the generator of $H^{2}(\SO(V);\Z) \cong \Z_{2}$, we have $[\tau] = -[\tau]$, so we can take $U$ to be a projective representation with cocycle $\tau$.  We now obtain
the Lie group/supergroup version of the Thom isomorphism:

\begin{proposition}\label{thom-group}
	The map $U \mapsto U \otimes \Cl(V)$, where $U$ is a $\tau$-twisted
	$G$-module and the action on the tensor product $U\otimes \Cl(V)$ is given by
	(\ref{eq:xi}), (\ref{eq:fi}), (\ref{eq:g}), gives an additive group isomorphism
	\begin{align*}
%	\RR{\tau}{}(G) \xrightarrow{\cong} \RR{b}{}(\Pi E), \qquad
	\SR{\tau}{*}(G) \xrightarrow{\cong} \SR{b}{*+\dim V}(\Pi E)
	\end{align*}
	from the twisted super representation group of $G$ to the $b$-twisted
	super representation group of $\Pi E$.
\end{proposition}

In terms of twisted equivariant $K$-theory (with compact supports), we have
$\SR{b}{*}(\Pi E) \cong K^{*}_{G}(V),$
and this isomorphism is precisely the equivariant Thom isomorphism $\sideset{^{\tau}}{^{*}_{G}}{\mathop{\!K}}(\mathrm{pt}) \cong K_{G}^{*+\dim V}(V)$ (see \cite{Ka}).
%Here, the Thom class is $[\Cl(V)] \in \SR{b+\tau}{\dim V}(\Pi E)$, and the Euler class is its restriction to $\g$:
%$$e = \bigl[\Cl_{0}(V)\bigr] - \bigl[\Cl_{1}(V)\bigr] = \bigl[\Lambda^{\mathrm{even}}(V \otimes\C)\bigr]
%- \bigl[\Lambda^{\mathrm{odd}}(V \otimes \C)\bigr] \, \in \, \SR{\tau}{}{}(G) \cong \K{\tau}{}{G}.$$
For even dimensional $V$, we can remove the degree shift in the Thom isomorphism by composing with the Bott periodicity isomorphism.
Doing so gives a Thom class $[\mathbb{S}^{V}]\in \SR{b+\tau}{}{}(\Pi E)$ and an Euler class $[\mathbb{S}^{V}_{0}] - [\mathbb{S}^{V}_{1}] \in \SR{\tau}{}{}(G) \cong \K{\tau}{}{G}$ in terms of the unique irreducible $\Cl(V)$-supermodule $\mathbb{S}^{V}= \mathbb{S}^{V}_{0}\oplus \mathbb{S}^{V}_{1}$. This Euler class plays a significant role in Sections \ref{Weyl-GKRS} and \ref{section-dirac-induction} below.
%(See Section \ref{pullbacks} below for a discussion of restrictions.)

For an alternative proof we can exponentiate the representation (\ref{eq:r}) to obtain an algebra homomorphism
$s : \dist_{\tau}(\widetilde{G}) \to \dist_{\tau}(\widetilde{G}) \otimes \ClR(V)$
in terms of the twisted convolution algebra  $\dist_{\tau}(\widetilde{G})$, induced by the group homomorphism
$s(g) = g \otimes \tilde{r}(g)$
for $g\in \widetilde{G}$.
This extends to an algebra isomorphism
$$f : \dist_{\tau}(\widetilde{G}) \semitensor \Cl(V) \xrightarrow{\cong} \dist_{\tau}(\widetilde{G}) \otimes \Cl(V),$$
given by
$$f : g \otimes \omega \mapsto r(g) (1 \otimes \omega) = g \otimes \tilde{r}(g) \cdot \omega$$
for $g\in \widetilde{G}$ and $\omega\in\Cl(V)$. Here, the multiplication on the semi-direct product $\dist_{\tau}(\widetilde{G}) \semitensor \ClR(V)$ is
$$(\phi\semitensor v)(\psi\semitensor w) = (\phi\conv\psi) \semitensor (v \cdot r_{\phi}w)
$$
for $\phi,\psi\in \dist_{\tau}(G)$ and $v,w\in\ClR(V)$,
which gives us our desired relation
$g \circ v \circ g^{-1} = r_{g} v$.

\section{Restriction and induction}
\label{restriction-induction}

\subsection{Pullbacks}
\label{pullbacks}
Let $f : A \to B$ be an (even) homomorphism of unital associative superalgebras.  If $r : B \to \End(V)$ is a representation of $B$ on a super vector space $V$, then
we can pull it back via $f$ to obtain a representation
$f^{*}r := r \circ f : A \to \End(V)$ of $A$ on $V$. If the $B$-action on
a $B$-module $V$ is understood, then we write $f^{*}V$ for the pulled back
$A$-module. Recalling our definition of the super representation group from Section \ref{subsection-representation-ring},
the pullback gives us a group homomorphism $f^{*} : F(B) \to F(A)$ on the free abelian groups generated by isomorphism classes of finite dimensional supermodules.  We also find that $f^{*}I(B) \subset I(A)$ for the subgroups
generated by classes $[U] - [V] + [W]$ whenever there exists a short exact
sequence
$$ 0 \to U \to V \to W \to 0$$
(with even maps).  Indeed, if we have such an exact sequence of super vector spaces which is $B$-equivariant, then it is clearly also equivariant with respect to the pulled back $A$-actions.  In addition, since $f: A\to B$ is an even map, the pullback of supermodules commutes with the parity reversal operator, $f^{*}\Pi = \Pi f^{*}$.  It follows that $f^{*}$ takes self-dual supermodules
$V \cong \Pi V$ to self-dual supermodules $f^{*}V \cong \Pi f^{*}V$,
%and anti-dual virtual supermodules $[V] - [\Pi V]$ to anti-dual supermodules $[f^{*}V] - [\Pi f^{*}V]$.
giving us %$f^{*}I_{-}(B) \subset I_{-}(A)$ and
$f^{*}I_{+}(B) \subset I_{+}(A)$, and thus $f^{*}$ descends to
a homomorphism
\begin{equation}\label{eq:pullback}
%	f^{*} :R(B) \to R(A), \qquad
	f^{*} :SR(B) \to SR(A)
\end{equation}
of the super representation groups.  Incorporating Clifford algebras, we can extend $f$ to a homomorphism
$f: A \tensor \Cl(n) \to B \tensor \Cl(n)$, and so we likewise obtain a pullback homomorphism
%\begin{equation}\label{eq:pullback-graded}
%	f^{*} : SR^{-n}(B) \to SR^{-n}(A)
%\end{equation}
of the degree-shifted super representation groups.  If $f : A \to B$ is a homomorphism of Hopf superalgebras, then the pullback (\ref{eq:pullback})
is a ring homomorphism, or a $\Z_{2}$-graded
ring homomorphism including degree shifts. If the homomorphism $f$ is an inclusion,
%be it of associative superalgebras, Lie groups, Lie superalgebras, Lie supergroups, or their twisted counterparts,
then we call the pullback the restriction.

A homomorphism $f: \h \to \g$ of Lie superalgebras extends via the injection $\g \to U(\g)$ to a Lie superalgebra homomorphism $f : \h \to U(\g)$ and then further extends to a homomorphism $f : U(\h) \to U(\g)$ of the universal enveloping algebras. Similarly, a homomorphism 
$f : H \to G$ of Lie groups extends to a homomorphism $f:\dist(H) \to \dist(G)$ of the convolution algebras (note that the convolution algebra is covariant as it is dual to the contravariant ring of smooth functions). The analogous statement holds for homomorphisms of Lie supergroups. We can therefore define pullback and restriction maps for the super representation rings of Lie superalgebras and Lie supergroups.

To construct pullbacks and restrictions of twisted super representation rings, we must pull back not only the representation but also the twisting.
For groups, recall from Section \ref{twist-group} that if $r : G \to \PU(V)$ is a projective representation of $G$, then its cocycle is the pullback
$\tau = r^{*}\bigl( \U(V) \bigr)$.
%\in H^{2}(G;\Z),$$ which further satisfies $d_{2}\tau = 0$, giving us an element $d_{3}\tau \in H^{3}(BG;\Z)$.
% then $r$ lifts to an true representation $\tilde{r} : \widetilde{G} \to \U(V)$, where $\widetilde{G} = r^{*}(\U(V))$ is classified by the cocycle $\tau = c_{1}(\widetilde{G}) = r^{*}c_{1}( \U(V) ) \in H^{2}(G)$, which satisfies $d_{2}\tau = 0$. Now,
If $f : H \to G$ is a Lie group homomorphism, then
$f^{*}r = r \circ f : H \to \PU(V)$ is a projective representation of $H$ with
cocycle
$(f^{*}r)^{*} \bigl(\U (V) \bigr) = f^{*} \bigl( r^{*}\bigl(\U(V)\bigr) \bigr) = f^{*} \tau.$
%Due to the naturality of the cohomology spectral sequence for a fibration, we have $d_{2}(f^{*}\tau) = 0$, and the corresponding class in $H^{3}(BH;\Z)$ is given by $d_{3}(f^{*}\tau) = f^{*}(d_{3}\tau)$.
%We therefore obtain a pullback homomorphism
%$$ f^{*} : \K{\tau}{}{G} \to \K{f^{*}\tau}{}{H}$$
%of the twisted Grothendieck groups.
%which when we take all twistings extends to a graded ring homomorphism  $f^{*}: \K{*}{}{G} \to \K{*}{}{H}$ with respect to the $H^{3}(BG;\Z)$-grading on $\K{*}{}{G}$ and the $H^{3}(BH;\Z)$-grading on $\K{*}{}{H}$.
For Lie superalgebras, a projective representation of $\g$ is a representation of a central extension $\tilde{\g}_{\tau}$ with cocycle $\tau$,
% represents a class $[\tau]\in H^{2}(\g)$,
such that the central generator $I$ acts by $i\,\mathrm{Id}$.  Given a Lie superalgebra homomorphism $f : \h \to \g$, it lifts to a Lie superalgebra homomorphism
$\tilde{f} : \tilde{\h}_{f^{*}\tau} \to \tilde{\g}_{\tau}$ which maps the central generator $I_{\tilde{h}}$ to $I_{\tilde{g}}$.  As a consequence, a $\tau$-twisted projective representation of $\g$ pulls back to a $f^{*}\tau$-twisted projective representation of $\h$.  In terms of the twisted universal enveloping algebra (\ref{eq:tuea}), we obtain an associative algebra homomorphism $f : U_{f^{*}\tau}(\h) \to U_{\tau}(\g)$. Likewise, a homomorphism $f: H\to G$ of Lie (super)groups induces an associative algebra homomorphism $f: \dist_{f^{*}\tau}(H) \to \dist_{\tau}(G)$ of the twisted convolution algebras defined in (\ref{eq:twisted-cg}). Such maps induce pullbacks or restrictions
\begin{alignat*}{3}
%f^{*} : \RR{\tau}{}(\g) &\to \RR{f^{*}\tau}{}(\h), &\qquad
  f^{*} : \SR{\tau}{*}(\g) &\to \SR{f^{*}\tau}{*}(\h) &\qquad
%f^{*} : \RR{\tau}{}(G) &\to \RR{f^{*}\tau}{}(H), &\qquad
  f^{*} : \SR{\tau}{*}(G) &\to \SR{f^{*}\tau}{*}(H).
\end{alignat*}
of the twisted super representation groups.

\subsection{Pushforwards}
\label{pushforwards}

%Recall our definition of the super representation group in Section \ref{subsection-representation-ring}, let us define for a unital associative algebra $A$ the Grothendieck group $R_{\Z_{2}}(A) := F(A) / I(A)$
%of (even) isomorphism classes of finite dimensional $A$-supermodules.  Since
%we are working with finite dimensional supermodules, the class of such a
%supermodule in $R_{\Z_{2}}(A)$ decomposes uniquely as the sum of classes
%of irreducible supermodules.  Thus $R_{\Z_{2}}(A)$ is the free abelian group generated by the classes of irreducibles.  Likewise, the representation group
%$R(A)$ and super representation group $SR(A)$ are free abelian groups on canonical bases of irreducibles described in Section \ref{subsection-representation-ring}.

If $A$ is a unital associative algebra and $V$ and $W$ are finite dimensional $A$-supermodules, we can define a symmetric bilinear pairing
%$$\langle V, W \rangle_{A} := \dim\Hom_A(V,W)_{0},$$
%which descends to $R_{\Z_{2}}(A)$ to give a symmetric bilinear form $R_{\Z_{2}}(A) \times R_{\Z_{2}}(A) \to \Z$.  Since we are counting the dimensions of only the \emph{even} homomorphisms, we see by Schur's Lemma that the classes of the irreducibles form an orthonormal basis for $R_{\Z_{2}}(A)$ with respect to this pairing. If we want to consider both even \emph{and odd} homomorphisms, we count the dimensions of the homomorphisms supersymmetrically, defining
$$\pair[V,W]_{A} := \sdim\Hom_{A}(V,W) = \dim\Hom_{A}(V,W)_{0} - \dim\Hom_{A}(V,W)_{1},$$
where $\sdim V = \dim V_{0} - \dim V_{1}$ is the superdimension,
in which we count both even \emph{and odd} homomorphisms supersymmetrically.
If $V \cong \Pi V$, then we see that $\pair[V,W]_{A} = 0$ for any $W$, and thus $\pair_{A}$ descends to a pairing on the super representation group $SR(A)$.  By Schur's Lemma, if $V$ and $W$ are irreducible $A$-supermodules, we have
$$\pair[V,W]_{A} = \begin{cases}
\phantom{-}1 & \text{ if $V,W$ are type M and $V \cong W$}, \\
-1 & \text{ if $V,W$ are type M and $V \cong \Pi W$}, \\
\phantom{-}0 & \text{ otherwise}.
\end{cases}$$
Recall that $SR(A)$ is a free abelian group on the (almost) canonical basis given in Section \ref{subsection-representation-ring}, consisting of one class $[M]$ for each pair $M,\Pi M$ of irreducibles of type M.  This basis is then orthonormal with respect to the
pairing $\pair_{A}$.\footnote{
Alternatively, the pairing 
$\langle\!\langle\!\langle\,{V,W}\,\rangle\!\rangle\!\rangle_{A} := \dim\Hom_{A}(V,W)_{0} + \dim\Hom_{A}(V,W)_{1}$
descends to the representation ring $R(A)$, but the canonical basis for $R(A)$ is \emph{not} orthonormal, as
$\langle\!\langle\!\langle \,Q,Q\,\rangle\!\rangle\!\rangle_{A} = 2$ for an irreducible of type Q.}

%\begin{remark}
%We can also define a pairing that descends to the representation group $R(A)$ by taking
%$$\langle\!\langle\!\langle\,{V,W}\,\rangle\!\rangle\!\rangle_{A} := \dim\Hom_{A}(V,W) = \dim\Hom_{A}(V,W)_{0} + \dim\Hom_{A}(V,W)_{1}.$$
%Here, the canonical basis for $R(A)$ described in Section \ref{subsection-representation-ring} is orthogonal, but \emph{not} orthonormal, as Schur's Lemma gives us $\langle\!\langle\!\langle \,Q,Q\,\rangle\!\rangle\!\rangle_{A} = 2$ for an irreducible $A$-supermodule $Q$ of type Q.  We will not consider this pairing further, and we will work exclusively with the super representation group.
%\end{remark}

In the spirit of Bott's paper \cite{B}, let $\widehat{SR}(A)$ denote the completion of $SR(A)$ with respect to the pairing $\pair_{A}$, by which we mean the additive group of formal, possibly infinite sums
\begin{equation}\label{eq:completion}
\widehat{SR}(A) := \Bigl\{ {\sum}_{i} a_{i} [V_{i}] \text{ for $a_{i}\in \Z$} \Bigr\},
\end{equation}
where $\{[V_{i}]\}$ is an orthonormal basis
for $SR(A)$ consisting of irreducibles.
The super representation group $SR(A)$ lies inside its completion $\widehat{SR}(A)$ as elements for which all but finitely many of the coefficients $a_{i}$ vanish.
Note that our pairing extends to a bilinear pairing $\pair_{A} : \widehat{SR}(A) \otimes SR(A) \to \Z$, but that we cannot in general pair two elements of the completion $\widehat{SR}(A)$. If $A$ is a Hopf algebra, then $SR(A)$ is a unital ring,
and the multiplication on $SR(A)$ extends to a bilinear map
$\widehat{SR}(A) \otimes SR(A) \to \widehat{SR}(A)$.  However, we do not in general
have a product on the completion $\widehat{SR}(A)$.  Classes in this completion $\widehat{SR}(A)$ can be represented by \emph{infinite} dimensional $A$-supermodules which satisfy the following finiteness conditions:
all their irreducible sub-supermodules must be finite dimensional, and each finite dimensional irreducible sub-supermodule appears with finite multiplicity.
%These are similar to finite type conditions, such as for positive energy representations of loop groups, in which an infinite dimensional module nevertheless decomposes into finite dimensional pieces.
%In this case, $\widehat{SR}(A)$ becomes a co-ring, with a comultiplication homomorphism
%$\widehat{SR}(A) \to \widehat{SR}(A) \otimes \widehat{SR}(A)$ (\cite{S}).

Given a superalgebra homomorphism $f: A \to B$, we construct a pushforward homomorphism $f_{*} : \widehat{SR}(A) \to \widehat{SR}(B)$ which satisfies the Frobenius reciprocity law:
\begin{equation}\label{eq:frob}
\pair[\,{[V]},f^{*}{[W]}\,]_{A} = \pair[\,f_{*}{[V]},{[W]}\,]_{B}
\end{equation}
for an $A$-supermodule $V$ and a $B$-supermodule $W$. In other words, we construct $f_{*}$ as the adjoint to the pullback $f^{*}$ with respect to the pairings $\pair$.  Given a class $[V] \in \widehat{SR}(A)$ in the completion,
in terms of a basis $\{[W_{i}]\}$ for $SR(B)$ of irreducibles, we define
\begin{equation}\label{eq:push-forward}
f_{*}[V] := {\sum}_{i} \pair[\,{[V]}, f^{*}{[W_{i}]}\,]_{A}\,[W_{i}].
\end{equation}
It follows immediately that this class $f_{*}[V]$ satisfies the Frobenius reciprocity law (\ref{eq:frob}). We note that this pushforward map is always an additive group homomorphism, \emph{not} a ring homomorphism.

\begin{remark}
Our homomorphism (\ref{eq:push-forward}) of super representation groups is based on the map used by Bott in \cite{B}, but it can also be realized explicitly in terms of representations as the pushforward $f_{*}V = B \otimes_{A} V$. Indeed, if $f : H \to G$ is a homomorphism of finite groups and $V$ is a finite dimensional representation of $H$, then its pushforward
$f_{*}V = \C[G] \otimes_{\C[H]} V$
is likewise finite dimensional. In such a case, the representation rings are themselves finite and thus remain unchanged upon completion. The two versions of the pushforward then agree: $[f_{*}V] = f_{*}[V]$. In general, the pushforward $f^{*}V$ is usually infinite dimensional, even if $V$ is finite dimensional, although these infinite dimensional representations do satisfy the finiteness conditions necessary for them to determine classes in the completed representation group. In this paper we prefer to work with finite dimensional representations, except to mention here the case underlying \cite{B}: If $i: H \hookrightarrow G$ is an inclusion of compact Lie groups, then the pushforward of a finite dimensional $H$-module $V$ is
$$i_{*}V = \dist(G) \otimes_{\dist(H)} V 
= \bigl( \dist(G) \otimes V \bigr)^{H} = \dist( G \times_{H} V ),$$
the space of distribution sections of the homogeneous vector bundle on $G/H$ induced by $V$. We will revisit this geometric interpretation in Section \ref{section-dirac-induction}.
%It follows from the Peter-Weyl theorem that this space of sections indeed satisfies the two finiteness conditions given above, and that its class $[f_{*}V]\in \widehat{SR}(G)$ in the completed super representation group agrees with our definition of $f_{*}[V]$ given in (\ref{eq:push-forward}). We will return to
\end{remark}

Introducing degree shifts via Clifford algebras, we obtain $\Z_{2}$-graded pushforward maps.
When working with twisted and degree-shifted representations of Lie superalgebras and Lie supergroups, we likewise obtain pushforward homomorphisms
$$f_{*} : \SRh{*}{*}(\h) \to \SRh{*}{*}(\g), \qquad
  f_{*} : \SRh{*}{*}(H) \to \SRh{*}{*}(G)$$
of additive groups, induced by Lie superalgebra homomorphisms $f: \h \to \g$ and Lie supergroup homomorphisms $f:H \to G$ respectively.
If the homomorphism $f$ is an inclusion,
%be it of associative superalgebras, Lie groups, Lie superalgebras, Lie supergroups, or their twisted counterparts, 
then we refer to the corresponding pushfoward map as the induction map.

\section{The Weyl-GKRS formula}
\label{Weyl-GKRS}

Let $\g$ be a semisimple Lie algebra, and let $\h$ be a reductive Lie subalgebra with $\mathrm{rank}\,\h = \mathrm{rank}\,\g$. Let $i : \h \hookrightarrow \g$ denote the inclusion of this equal rank subalgebra.
Let $b$ be an $\ad$-invariant inner product on $\g$, such as the Killing form,
with respect to which $\g$ decomposes as the direct sum of $\h$ and its orthogonal complement. The inner product $b$ gives
an $\ad^{*}_{\g}$-invariant inner product on the dual space $\g^{*}$, with
respect to which $\g^{*}$ decomposes as the direct sum of $\h^{*}$ and the dual
to the orthogonal complement of $\h$ in $\g$.  The inner product $b$ then restricts to an $\ad^{*}_{\h}$-invariant inner product on $\h^{*}$, which we also denote by $b$.
%With respect to $b$, we have an inclusion $\h^{*}\subset\g^{*}$.
Now, consider the coadjoint Lie superalgebra $\g \oplus \Pi \g^{*}$ and its
Lie sub-superalgebra $\h \oplus \Pi \h^{*}$, with the inclusion denoted by
$j : \h \oplus \Pi \h^{*} \to \g \oplus \Pi \g^{*}$.  

We can now consider the following diagram of maps:
\begin{equation}\label{eq:noncommute-algebra}
\begin{CD}
	SR(\g) @>\mathrm{Thom}_{\g}>> \SR{b}{\dim\g}(\g\oplus\Pi \g^{*}) \\
	@Vi^{*}VV @VVj^{*}V \\
	SR(\h)  @>\mathrm{Thom}_{\h}>> \SR{b}{\dim\h}(\h\oplus\Pi \h^{*})
\end{CD}
%\begin{diagram}
%	SR(\g) & \rTo^{\mathrm{Thom}_{\g}} & \SR{b}{\dim\g}(\g\oplus\Pi \g^{*}) \\
%	\dTo<{i^{*}} & & \dTo>{j^{*}} \\
%	SR(\h) & \rTo^{\mathrm{Thom}_{\h}} & \SR{b}{\dim\h}(\h\oplus\Pi \h^{*})
%\end{diagram}
\end{equation}
where the horizontal maps are Thom isomorphisms as described in Section \ref{thom-algebra}, and the vertical maps are restriction maps as described
in Section \ref{pullbacks}.  Note that since $\h$ has maximal rank in $\g$,
we have $\dim \h \equiv \dim \g \pmod 2$, so the degree shifts on the right hand
side of the diagram agree. (Strictly speaking, the map on the right side
is $j^{*}$ composed with the Bott periodicity isomorphism $\SR{b}{\dim\g}(\g\oplus\Pi\g^{*}) \to \SR{b}{\dim\h}(\g\oplus\Pi\g^{*})$.)
This diagram does \emph{not} commute!  In fact, we have

\begin{theorem}\label{theorem-restriction}
The Lie superalgebra restriction map $j^{*} : \SR{b}{\dim\g}(\g\oplus\Pi \g^{*}) \to \SR{b}{\dim\h}(\h\oplus\Pi \h^{*})$ pulls
back via the Thom isomorphisms to a map $SR(\g) \to SR(\h)$,
$$SR(\g) \xrightarrow{\mathrm{Thom}_{\g}} \SR{b}{\dim\g}(\g\oplus\Pi\g^{*})
\xrightarrow{j^{*}\,\circ\,\mathrm{Bott}} \SR{b}{\dim\h}(\h\oplus\Pi\h^{*})
\xrightarrow{(\mathrm{Thom}_{\h})^{-1}} SR(\h),$$
which is given by
\begin{equation}\label{eq:gkrs-map}
[V] \in SR(\g) \longmapsto i^{*}[V] \, \bigl( [\mathbb{S}_{0}] - [\mathbb{S}_{1}] \bigr) \in SR(\h),
\end{equation}
where $\mathbb{S} = \mathbb{S}_{0}\oplus \mathbb{S}_{1}$ is the unique irreducible $\Cl(\g^{*}/\h^{*},b)$-supermodule.
\end{theorem}

\begin{proof}
We begin with a $\g$-supermodule $V$ corresponding to $[V]\in SR(\g)$.
Under the Thom isomorphism, it maps to the $b$-twisted $\g\oplus\Pi\g^{*}$-Clifford supermodule $\mathrm{Thom}_{\g}\,V = V \otimes \Cl(\g^{*},b)$, with actions
(\ref{eq:X}), (\ref{eq:xi}), (\ref{eq:fi}). The Clifford algebra is multiplicative, factoring $\h$-equivariantly as
\begin{equation}\label{eq:clifford-split}
	\Cl(\g^{*},b) \cong \Cl(\h^{*},b) \tensor \Cl(\g^{*}/\h^{*},b)
	\cong \Cl(\h^{*},b) \otimes \Cl(\g^{*}/\h^{*},b),
\end{equation}
where for the second isomorphism we use the twofold periodicity of complex Clifford algebras to change the graded tensor product into an ungraded one.
Since $\g^{*}/\h^{*}$ is even dimensional, its complex Clifford algebra is $\Cl(\g^{*}/\h^{*},b) \cong \End(\mathbb{S})$, where $\mathbb{S} = \mathbb{S}_{0}\oplus\mathbb{S}_{1}$
is the unique irreducible complex Clifford supermodule up to isomorphism (see \cite{LM}). The space of endomorphisms then decomposes as the tensor product $\End(\mathbb{S}) \cong \mathbb{S} \otimes \mathbb{S}^{*}$ with respect to the left and right
actions. Combining this with (\ref{eq:clifford-split}), we obtain
$$\mathrm{Thom}_{\g}\,V \cong V \otimes \Cl(\h^{*},b)\otimes \mathbb{S} \otimes \mathbb{S}^{*},$$
where $\g\oplus\Pi \g^{*}$ and $\Cl(\dim\h)$ act on the $(V \otimes \Cl(\h^{*}) \otimes \mathbb{S})$ factors, while only $\Cl(\dim\g - \dim\h)$ acts on the $\mathbb{S}^{*}$ factor.  Applying the Bott periodicity isomorphism
eliminates the $\mathbb{S}^{*}$ factor, and then restricting to 
$\h \oplus \Pi \h^{*}$, we regroup the terms to obtain:
$$j^{*} \circ \mathrm{Bott} \circ \mathrm{Thom}_{\g} (V) 
\cong ( i^{*}V \otimes \mathbb{S} ) \otimes \Cl(\h^{*},b).$$
Undoing the last Thom isomorphism gives us
$$(\mathrm{Thom}_{\h})^{-1} \circ j^{*} \circ \mathrm{Bott} \circ \mathrm{Thom}_{\g} (V) \cong i^{*}V \otimes \mathbb{S},$$
and finally, since $i^{*}V \otimes \mathbb{S} = 
(i^{*}V \otimes \mathbb{S}_{0}) \oplus (i^{*}V\otimes \mathbb{S}_{1})$ is an
$\h$-supermodule, its class in the super representation ring is that of
the virtual $\h$-module
$$[i^{*}V \otimes \mathbb{S}_{0}] - [i^{*}V \otimes \mathbb{S}_{1}]
 = i^{*}[V] \bigl( [\mathbb{S}_{0}] - [\mathbb{S}_{1}] \bigr) \in SR(\h), $$
given by taking the formal direct difference of its even and odd components as in (\ref{eq:direct-difference}).
%as we found in our discussion of the purely even case in Section \ref{even}.
\end{proof}

If $V$ is an irreducible $\g$-module, then its image under the map (\ref{eq:gkrs-map}) was computed explicitly by Gross, Kostant, Ramond, and Sternberg in \cite{GKRS}. Their formula, given in the following Theorem, reduces to the Weyl character formula when $\h = \mathfrak{t}$ is a Cartan subalgebra. Here, we assume that we have chosen a common Cartan subalgebra $\mathfrak{t}\subset\h\subset\g$, giving us a compatible Weyl groups $W_{\h}\subset W_{\g}$, and that we have chosen compatible systems of positive roots for $\h$ and $\g$.

\begin{theorem}[GKRS]\label{GKRS}
If $V_{\lambda}$ is an irreducible $\g$-module with highest weight $\lambda$,
then
$$ [i^{*}V_{\lambda} \otimes \mathbb{S}_{0}] - [i^{*}V_{\lambda} \otimes \mathbb{S}_{1}]
   = {\sum}_{c\in W_{\g}/W_{\h}} (-1)^{c}\,\bigl[ U_{c(\lambda + \rho_{\g}) - \rho_{\h}} \bigr] \in \K{}{}{\g},$$
where $U_{\mu}$ denotes the irreducible $\h$-module with highest weight $\mu$,
the weights $\rho_{\g}$ and $\rho_{\h}$ are half the sum of the positive roots
of $\g$ and $\h$ respectively, and the representative $c$ of each coset in $W_{\g}/W_{\h}$ is chosen so that the weight $c(\lambda+\rho_{\g})-\rho_{\h}$ on the right hand side is dominant.
\end{theorem}

We also have a Lie group version of Theorem \ref{theorem-restriction}. Let $G$ be a (connected) compact semisimple Lie group with Lie algebra $\g$, and let $H$, with Lie algebra $\h$, be an equal rank Lie subgroup of $G$.  Let $i: H\hookrightarrow G$ denote the inclusion.  We can then construct their parity reversed cotangent bundles, the Lie supergroups $\Pi(T^{*}G)$ and $\Pi(T^{*}H)$ with underlying even Lie groups $G$ and $H$ and associated Lie superalgebras $\g\oplus\Pi\g^{*}$ and $\h\oplus\Pi\h^{*}$ respectively. Let $j : \Pi(T^{*}H) \to \Pi(T^{*}G)$ denote the inclusion of the Lie supergroups (actually corresponding to inclusions of their underlying even Lie groups and associated Lie superalgebras). As before,
we choose an $\Ad^{*}_{G}$-invariant inner product $b$ on $\g^{*}$, which restricts to an $\Ad^{*}_{H}$-invariant inner product on $\h^{*}$ which we
also denote by $b$.  This gives us our twistings on the odd components.
We
also have twistings $\tau_{G} = \Ad^{*}\Spinc(\g^{*})$
and $\tau_{H} = \Ad^{*}\Spinc(\h^{*})$ on the even components, constructed as in Section \ref{section-thom-group}.
%given by (\ref{eq:spinc-twisting}), taking $V = \g^{*}$ and $V = \h^{*}$ to be the coadjoint representations of $G$ and $H$ respectively.
%We
%also have twistings $\tau_{G}\in H^{2}(G;\Z)$
%and $\tau_{H}\in H^{2}(H;\Z)$ (satisfying $d_{2}\tau_{G} = d_{2}\tau_{H} = 0$) on the even components, given by (\ref{eq:spinc-twisting}), taking $V = \g^{*}$
%and $V = \h^{*}$ to be the coadjoint representations of $G$ and $H$ respectively.
%, as well as a possible auxiliary twistings given by $\sigma\in H^{2}(G;\Z)$ with $d_{2}\sigma = 0$
The Lie group counterpart of (\ref{eq:noncommute-algebra}) is the non-commutative diagram:
\begin{equation*}
\begin{CD}
	\SR{\tau_{G}}{}(G) @>\mathrm{Thom}_{G}>> \SR{b}{\dim G}\bigl(\Pi(T^{*}G)\bigr) \\
	@Vi^{*}VV @VVj^{*}V \\
	\SR{\tau_{H}}{}(H) @>\mathrm{Thom}_{H}>> \SR{b}{\dim H}\bigl(\Pi(T^{*}H)\bigr)
\end{CD}
%\begin{diagram}
%	\SR{\tau_{G}}{}(G) & \rTo^{\mathrm{Thom}_{G}} & \SR{b}{\dim G}\bigl(\Pi(T^{*}G)\bigr) \\
%	\dTo<{i^{*}} & & \dTo>{j^{*}} \\
%	\SR{\tau_{H}}{}(H) & \rTo^{\mathrm{Thom}_{H}} & \SR{b}{\dim H}\bigl(\Pi(T^{*}H)\bigr)
%\end{diagram}
\end{equation*}
and the Lie group counterpart of Theorem \ref{theorem-restriction} is then:

\begin{theorem}\label{supergroup-restriction}
The Lie supergroup restriction map $j^{*} : \SR{b}{\dim G}\bigl(\Pi(T^{*}G)\bigr) \to \SR{b}{\dim H}\bigl(\Pi(T^{*}H)\bigr)$ pulls
back via the Thom isomorphisms to a map $\SR{\tau_{G}}{}(G) \to \SR{\tau_{H}}{}(H)$,
$$\SR{\tau_{G}}{}(G) \xrightarrow{\mathrm{Thom}_{G}} \SR{b}{\dim G}\bigl(\Pi(T^{*}G)\bigr)
\xrightarrow{j^{*}\,\circ\,\mathrm{Bott}} \SR{b}{\dim H}\bigl(\Pi(T^{*}H)\bigr)
\xrightarrow{(\mathrm{Thom}_{H})^{-1}} \SR{\tau_{H}}{}(H),$$
which is given by
\begin{equation*}
[V] \in \SR{\tau_{G}}{}(G) \longmapsto i^{*}[V] \, \bigl( [\mathbb{S}_{0}] - [\mathbb{S}_{1}] \bigr) \in \SR{\tau_{H}}{}(H),
\end{equation*}
where $\mathbb{S} = \mathbb{S}_{0}\oplus \mathbb{S}_{1}$ is the unique irreducible $\Cl(\g^{*}/\h^{*},b)$-supermodule.
\end{theorem}

\begin{remark}
We note that the class of the cocycle $[\tau_{G}]\in H^{2}(G;\Z)$ satisfies $2[\tau_{G}]= 0$. So, if $G$ is simply connected, or if $\pi_{1}(G)$ has no $2$-torsion, then the class $[\tau_{G}]$ must vanish, and we can deal with actual $G$-modules rather than projective ones.  Also, if $H=T$ is a maximal torus, then $H^{3}(BT;
\Z) = 0$, and thus the class $[\tau_{T}]\in H^{2}(T;\Z)$ vanishes. 
\end{remark}

\section{Dirac induction}
\label{section-dirac-induction}

We continue to use the notation of the previous section. Replacing the restriction maps in (\ref{eq:noncommute-algebra}) with induction maps, we obtain a non-commuting diagram
$$\begin{CD}
	\SRh{}{}(\g) @>\mathrm{Thom}_{\g}>> \SRh{b}{\dim\g}(\g\oplus\Pi \g^{*}) \\
	@Ai_{*}AA @AAj_{*}A \\
	\SRh{}{}(\h) @>\mathrm{Thom}_{\h}>>  \SRh{b}{\dim\h}(\h\oplus\Pi \h^{*})
\end{CD}
$$
%\begin{diagram}
%	\SRh{}{}(\g) & \rTo^{\mathrm{Thom}_{\g}} & \SRh{b}{\dim\g}(\g\oplus\Pi \g^{*}) \\
%	\uTo<{i_{*}} & & \uTo>{j_{*}} \\
%	\SRh{}{}(\h) & \rTo^{\mathrm{Thom}_{\h}} & \SRh{b}{\dim\h}(\h\oplus\Pi \h^{*})
%\end{diagram}
where the Thom isomorphisms along the top and bottom are the extensions of the Thom isomorphisms on the super representation rings $\SR{}{}$ to their completions $\SRh{}{}$.

\begin{theorem}\label{theorem-induction}
The Lie superalgebra induction map $j_{*} : \SRh{b}{\dim\h}(\h\oplus\Pi \h^{*}) \to \SRh{b}{\dim\g}(\g\oplus\Pi \g^{*})$ pulls
back via the Thom isomorphisms to a map $\SRh{}{}(\h)\to \SRh{}{}(\g)$,
$$\SRh{}{}(\h) \xrightarrow{\mathrm{Thom}_{\h}} \SRh{b}{\dim\h}(\h\oplus\Pi\h^{*})
\xrightarrow{j_{*}\,\circ\,\mathrm{Bott}} \SRh{b}{\dim\g}(\g\oplus\Pi\g^{*})
\xrightarrow{(\mathrm{Thom}_{\g})^{-1}} \SRh{}{}(\g),$$
which is given by
\begin{equation}\label{eq:dirac-map}
[U] \in \SRh{}{}(\h) \longmapsto i_{*} [U \otimes \mathbb{S}_{0}^{*}] - i_{*}[U \otimes \mathbb{S}_{1}^{*}] \in \SRh{}{}(\g),
\end{equation}
where $\mathbb{S} = \mathbb{S}_{0}\oplus \mathbb{S}_{1}$ is the unique irreducible $\Cl(\g^{*}/\h^{*},b)$-supermodule.
\end{theorem}

\begin{proof}
Consider classes $[U] \in \SRh{}{}(\h)$ and $[V]\in\SR{}{}(\g)$.  Applying the Thom isomorphisms, we obtain classes $\mathrm{Thom}_{\h}[U]\in\SRh{}{}(\h\oplus\Pi\h^{*})$ and
$\mathrm{Thom}_{\g}[V]\in\SR{}{}(\g\oplus\Pi\g^{*})$.  Using the pairing
$\SRh{}{} \otimes \SR{}{} \to \Z$, the Frobenius reciprocity law (\ref{eq:frob})
gives us
\begin{equation*}
\pair[j_{*}\,\mathrm{Thom}_{\h}\,{[U]}, \, \mathrm{Thom}_{\g}\, {[V]}]_{\g\oplus\Pi\g^{*}}
=
\pair[\mathrm{Thom}_{\h}\, {[U]}, \, j^{*}\,\mathrm{Thom}_{\g}\,{[V]}]_{\h\oplus\Pi\h^{*}}
.
\end{equation*}
Since the Thom isomorphisms preserve the pairing, we can perform this
computation on the Lie algebras $\g$ and $\h$ instead of their Lie superalgebra counterparts, and we obtain
\begin{equation*}
	\pair[ (\mathrm{Thom}_{\g})^{-1} j_{*}\,\mathrm{Thom}_{\h}\,{[U]}, \,
	{[V]} ]_{\g}
	=
	\pair[ {[U]}, \, (\mathrm{Thom}_{\h})^{-1} j^{*}\,\mathrm{Thom}_{\g}\, 
	{[V]} ]_{\h}
	.
\end{equation*}
Finally, applying Theorem \ref{theorem-restriction}, we have
\begin{equation*}\begin{split}
	\pair[ (\mathrm{Thom}_{\g})^{-1} j_{*}\,\mathrm{Thom}_{\h}\,{[U]}, \,
	{[V]} ]_{\g}
	&= \pair[ {[U]}, \, i^{*} {[V]}\,\bigl( {[\mathbb{S}_{0}]} - {[\mathbb{S}_{1}]}\bigr)]_{\h} \\
	&= \pair[ {[U]}\,\bigl( {[\mathbb{S}_{0}^{*}]} - {[\mathbb{S}_{1}^{*}]}\bigr), \, i^{*} {[V]}]_{\h} \\
	&= \pair[ i_{*} \bigl( {[U \otimes \mathbb{S}_{0}^{*}]} - {[U \otimes \mathbb{S}_{1}^{*}]}\bigr), \, {[V]}]_{\g}.
\end{split}\end{equation*}
Since $(\mathrm{Thom}_{\g})^{-1}j_{*}\mathrm{Thom}_{\h}[U] \in \SRh{}{}(\g)$,
it is completely determined by its pairings with the classes $[V_{i}]$ of irreducibles in $SR(\g)$, and thus $(\mathrm{Thom}_{\g})^{-1}j_{*}\mathrm{Thom}_{\h}[U] = i_{*}[U\otimes \mathbb{S}_{0}^{*}] - i^{*}[U \otimes S_{1}^{*}]$.
\end{proof}

We also have a version for Lie groups, which incorporates the twistings $\tau_{G}$ and $\tau_{H}$ discussed in the last section.  In the Lie group case
we consider the non-commuting diagram:
\begin{equation*}
\begin{CD}
	\SRh{\tau_{G}}{}(G)@>\mathrm{Thom}_{G}>> \SRh{b}{\dim G}\bigl(\Pi(T^{*}G)\bigr) \\
	@Ai_{*}AA @AAj_{*}A \\
	\SRh{\tau_{H}}{}(H) @>\mathrm{Thom}_{H}>> \SRh{b}{\dim H}\bigl(\Pi(T^{*}H)\bigr)
\end{CD}
%\begin{diagram}
%	\SR{\tau_{G}}{}(G) & \rTo^{\mathrm{Thom}_{G}} & \SR{b}{\dim G}\bigl(\Pi(T^{*}G)\bigr) \\
%	\uTo<{i_{*}} & & \uTo>{j_{*}} \\
%	\SR{\tau_{H}}{}(H) & \rTo^{\mathrm{Thom}_{H}} & \SR{b}{\dim H}\bigl(\Pi(T^{*}H)\bigr)
%\end{diagram}
\end{equation*}
and our theorem takes the following form:

\begin{theorem}\label{group-induction}
The Lie supergroup induction map $j_{*} : \SRh{b}{\dim H}\bigl(\Pi(T^{*}H)\bigr) \to \SRh{b}{\dim G}\bigl(\Pi(T^{*}G)\bigr)$
pulls back via the Thom isomorphisms to a map $\SRh{\tau_{H}}{}(H) \to \SRh{\tau_{G}}{}(G)$,
$$\SRh{\tau_{H}}{}(H) \xrightarrow{\mathrm{Thom}_{H}} \SRh{b}{\dim H}\bigl(\Pi(T^{*}H)\bigr)
\xrightarrow{j^{*}\,\circ\,\mathrm{Bott}} \SRh{b}{\dim G}\bigl(\Pi(T^{*}G)\bigr)
\xrightarrow{(\mathrm{Thom}_{G})^{-1}} \SRh{\tau_{G}}{}(G)$$
which is given by
\begin{equation*}
[U] \in \SRh{\tau_{H}}{}(H) \longmapsto i_{*} [U \otimes \mathbb{S}_{0}^{*}] - i_{*}[U \otimes \mathbb{S}_{1}^{*}] \in \SRh{\tau_{G}}{}(G),
\end{equation*}
where $\mathbb{S} = \mathbb{S}_{0}\oplus \mathbb{S}_{1}$ is the unique irreducible $\Cl(\g^{*}/\h^{*},b)$-supermodule.
\end{theorem}

In \cite{B}, Bott defines the induction map $i_{*}$ for Lie group representations, giving a geometric interpretation in terms of the index of homogeneous elliptic operators.  Recall the isomorphism $R(H) \to K_{G}( G/H )$, which assigns to a finite dimensional $H$-module $U$ the homogeneous vector bundle $G\times_{H}U$ over
the coset space $G/H$.  The space of $L^{2}$-sections $\Gamma(G \times_{H}U)$
is an infinite dimensional $G$-module satisfying the finiteness conditions of Section \ref{pushforwards}, whose class in the completed representation ring gives the induced representation:
$$\bigl[ \Gamma(G \times_{H}U)\bigr] = i_{*}[H] \in \SRh{}{}(G).$$
Indeed, for each finite dimensional $G$-module $V$, we have the Frobenius reciprocity isomorphism
\begin{equation}\label{eq:frob1}
\Hom_{G}\bigl(\,\Gamma(G\times_{H}U), \, V\,\bigr) \cong
  \Hom_{H}( U, i^{*}V ),
\end{equation}
corresponding to the defining identity (\ref{eq:frob}) for the induction map. In addition, it follows from the Peter-Weyl theorem that (\ref{eq:frob1}) gives a
complete description of $\Gamma(G\times_{H}U)$. In particular, all of its irreducible subrepresentations are finite dimensional.

Bott's theorem then says that the index of an elliptic homogeneous differential operator is simply the difference of its domain and codomain, regardless of the operator itself:

\begin{theorem}[Bott]\label{bott}
Given $H$-modules $U_{1}$ and $U_{2}$, if $D : \Gamma(G\times_{H} U_{1}) \to \Gamma(G \times_{H}U_{2})$ is an elliptic homogeneous differential operator, then its $G$-equivariant index is 
$$\mathrm{Index}_{G}\, D
= \bigl[ \Gamma(G\times_{H}U_{1}) \bigr]
- \bigl[ \Gamma(G\times_{H}U_{2}) \bigr]
= i_{*}\bigl([U_{1}]-[U_{2}]\bigr) \in \SRh{}{}(G),$$
and furthermore, the index is actually a finite element in $\SR{}{}(G) \subset \SRh{}{}(G)$.
\end{theorem}

In particular, if $G/H$ is spin, we can consider the Dirac operator $\dirac^{G/H}$.  The tangent bundle of $G/H$ is $T(G/H) \cong G \times_{H} (\g/\h)$, and it follows
that the spin bundle is $S \cong G \times_{H} \mathbb{S}$, which decomposes into the two half-spin bundles $S =  S^{+}\oplus S^{-}$ given by
$S^{+} \cong G \times_{H} \mathbb{S}_{0}$ and
$S^{-} \cong G \times_{H} \mathbb{S}_{1}$.  Given any finite dimensional $H$-module $U$, we can consider the Dirac operator with values in $G\times_{H}U$,
\begin{equation}\label{eq:homogeneous-dirac}
\dirac^{G/H}_{U} : \Gamma \bigl( G \times_{H} ( \mathbb{S}_{0}^{*}\otimes U ) \bigr) \to
\Gamma \bigl( G \times_{H} ( \mathbb{S}_{1}^{*}\otimes U ) \bigr).
\end{equation}
This is an elliptic homogeneous differential operator, so applying Bott's Theorem we obtain
\begin{equation}\label{eq:index-formula}
\mathrm{Index}_{G}\,\dirac^{G/H}_{U} = i_{*}[U \otimes \mathbb{S}_{0}^{*}] -
i_{*}[U \otimes \mathbb{S}_{1}^{*}] \in \SR{}{}(G).
\end{equation}
Even if $G/H$ is not spin, when it is $\mathrm{Spin}^{c}$
we can construct Dirac operators of the form (\ref{eq:homogeneous-dirac})
for $(\tau_{H}-i^{*}\tau_{G})$-twisted projective representations.
%In this case, the $G$-index (\ref{eq:index-formula}) takes values in the representation ring $\SR{\tau_{G}}{}(G)$.
Comparing the index (\ref{eq:index-formula}) with Theorem \ref{group-induction}, we get
%(after shifting the twistings by $\tau_{G}$):

\begin{theorem}\label{dirac-induction}
The Lie supergroup induction map restricts to an additive homomorphism
on the uncompleted representation groups,
$j_{*} : \SR{b}{\dim H}\bigl(\Pi(T^{*}H)\bigr) \to \SR{b}{\dim G}\bigl(\Pi(T^{*}G)\bigr)$, which pulls back via the Thom isomorphisms to$$(\mathrm{Thom}_{G})^{-1} \circ j_{*} \circ \mathrm{Thom}_{H} = \mathrm{Index}_{G}\,\dirac^{G/H}_{\bullet} : \SR{\tau_{H}-i^{*}\tau_{G}}{}(H) \to \SR{}{}(G),$$
the Dirac induction map.
\end{theorem}

This Dirac induction map is related to the holomorphic induction map of the Borel-Weil-Bott theorem, which takes an $H$-module $U$ to the space of holomorphic sections of $G \times_{H} U$. In order to apply holomorphic induction, the coset space $G/H$ must be a complex homogeneous space, such as occurs when $H=T$ is a maximal torus in $G$. In contrast, Dirac induction requires only a $\mathrm{Spin}^{c}$-structure on $G/H$, which is a weaker condition. On the other hand, the space of holomorphic sections is just the degree 0 component of the Dolbeault cohomology $H_{\bar{\partial}}^{*}(G/H;U)$, which associates a $\Z$-graded $G$-module to each $H$-module $U$. When considering harmonic spinors, we obtain only a $\Z_{2}$-graded $G$-module consisting of the kernel and cokernel of the Dirac operator $\dirac^{G/H}_{U}$.  The $G$-index (\ref{eq:index-formula}) then recovers the Euler characteristic of the Dolbeault complex (up to a $\rho$-shift resulting from tensoring with the canonical complex line bundle of $G/H$).  Thus, the presence of a complex structure on $G/H$ allows us to extend the $\Z_{2}$-grading to a richer $\Z$-grading.

In analogy to the Borel-Weil-Bott theorem, we can use Theorem \ref{GKRS} to compute the Dirac induction map applied to the class of an irreducible $H$-module $U_{\mu}$ with highest weight $\lambda$.  Using the same notation as we developed for Theorem \ref{GKRS}, we obtain:

\begin{theorem}\label{dirac-index}
	If $U_{\mu}$ is an irreducible $H$-module with highest weight $\mu$,
	then
	$$(\mathrm{Thom}_{G})^{-1}\, j_{*}\, \mathrm{Thom}_{H}[U_{\mu}]
	= \mathrm{Index}_{G}\,\dirac^{G/H}_{U_{\mu}} = (-1)^{c}\, \bigl[V_{c(\mu+\rho_{H}) - \rho_{G}}\bigr],$$
	if there exists a Weyl group element $c\in W_{G}$ such that $c(\mu+\rho_{H}) - \rho_{G}$ is dominant, or $0$ otherwise.
\end{theorem}

See \cite{L0} for a quick proof  based on Theorem \ref{GKRS}. This result is a weak version of the Borel-Weil-Bott theorem which holds
for all equal rank subgroups $H\subset G$, not just the cases where $G/H$ is complex.  While the Borel-Weil-Bott theorem gives us an induced representation in a specific integer degree in the Dolbeault cohomology, the representation given by Dirac induction carries only a $\Z_{2}$-degree given by the sign $(-1)^{c}$. Nevertheless, this Dirac induction is extremely useful in representation theory, as it allows us to explicitly construct any finite dimensional $G$-module as a space of harmonic spinors on $G/H$.

\begin{remark}
	In \cite{Sl1,Sl2}, Slebarski proved a stronger version of Theorem \ref{dirac-index}
	by computing not just the index, but in fact the kernel and cokernel
	of the Dirac operator for a 1-parameter family of connections on $G/H$.
	For one particular choice of connection, referred to as the
	``reductive connection'' by Slebarski and constructed independently
	as part of the ``cubic'' Dirac operator by Alekseev and Meinrenken in \cite{AM} and Kostant in \cite{K0}, we find that the index lies completely
	in the kernel or the cokernel, giving a vanishing theorem similar to
	that of the Borel-Weil-Bott theorem.  Indeed, Kostant uses his cubic
	Dirac operator in \cite{K} to give an alternative proof of the Lie algebra 
	cohomology version of the Borel-Weil-Bott theorem.
\end{remark}

\end{document}